\documentclass[11pt]{article} % use larger type; default would be 10pt

\usepackage[utf8]{inputenc} % set input encoding (not needed with XeLaTeX)

\usepackage{geometry} % to change the page dimensions
\usepackage{graphicx} % support the \includegraphics command and options

\usepackage{booktabs} % for much better looking tables
\usepackage{array} % for better arrays (eg matrices) in maths
\usepackage{paralist} % very flexible & customisable lists (eg. enumerate/itemize, etc.)
\usepackage{verbatim} % adds environment for commenting out blocks of text & for better verbatim
\usepackage{subfig} % make it possible to include more than one captioned figure/table in a single float
\usepackage{fancyhdr} % This should be set AFTER setting up the page geometry
\usepackage{sectsty}
\allsectionsfont{\sffamily\mdseries\upshape} % (See the fntguide.pdf for font help)
\usepackage[nottoc,notlof,notlot]{tocbibind} % Put the bibliography in the ToC
\usepackage[titles,subfigure]{tocloft} % Alter the style of the Table of Contents

\usepackage{amsmath, amssymb, bm}
\usepackage[T1]{fontenc}
\usepackage{textcomp}

\newtheorem{theorem}{Theorem}[section]
\newtheorem{lemma}[theorem]{Lemma}
\newtheorem{proposition}[theorem]{Proposition}

\newtheorem{definition}[theorem]{Definition}

\numberwithin{equation}{section} 

\title{Effect of microscopic pausing time distributions \\ 
on the dynamical limit shapes for random Young diagrams}
\author{Akihito HORA}
\date{Department of Mathematics, Faculty of Science, Hokkaido University} 
\begin{document}
\maketitle

\begin{abstract}
The irreducible decomposition of successive restriction and induction of irreducible representations of a symmetric group 
gives rise to a Markov chain on Young diagrams keeping the Plancherel measure invariant. 
Starting from this Res-Ind chain, we introduce a not necessarily Markovian continuous time random walk on Young diagrams 
by considering a general pausing time distribution between jumps according to the transition probability 
of the Res-Ind chain. 
We show that, under appropriate assumptions for the pausing time distribution, a diffusive scaling limit brings us 
concentration at a certain limit shape depending on macroscopic time which leads to a similar consequence 
to the exponentially distributed case studied in our earlier work. 
The time evolution of the limit shape is well described by using free probability theory. 
On the other hand, we illustrate an anomalous phenomenon observed with a pausing time obeying 
a one-sided stable distribution, heavy-tailed without the mean, 
in which a nontrivial behavior appears under a non-diffusive regime of the scaling limit.
\end{abstract}

\section{Introduction}

As a remarkable classical result in the field of asymptotic representation theory, 
the limit shape of random Young diagrams originated with Vershik--Kerov \cite{VeKe77} and Logan--Shepp \cite{LoSh77}. 
Let $\mathbb{Y}$ denote the set of Young diagrams. 
Set $\mathbb{Y}_n = \bigl\{\lambda\in\mathbb{Y} \big| |\lambda| =n\bigr\}$, 
where $|\lambda|$ denotes the size ($=$ the number of boxes) of $\lambda\in\mathbb{Y}$. 
For $\lambda = (\lambda_1\geqq \lambda_2 \geqq \cdots) \in \mathbb{Y}$, set 
$m_j(\lambda) = \sharp \{ i | \lambda_i = j\}$, namely the number of rows of length $j$ in $\lambda$. 
The total number of rows is $l(\lambda) = \sum_{j=1}^\infty m_j(\lambda)$. 
The Plancherel measure on $\mathbb{Y}_n$ is defined by 
%
%\begin{equation}\label{eq:1-1}
\[ 
\mathrm{M}_{\mathrm{Pl}}^{(n)} (\{\lambda\}) = \frac{(\dim\lambda)^2}{n!}, \qquad \lambda\in\mathbb{Y}_n.
\] 
%\end{equation}
%
Young diagram $\lambda$ is identified with the profile $y = \lambda(x)$ depicted in the $xy$ coordinates plane, satisfying 
\[
\int_{\mathbb{R}}\bigl( \lambda(x) - |x|\bigr) dx = 2 |\lambda|
\] 
(see Appendix, Figure~\ref{fig:2}). 
For $\lambda\in\mathbb{Y}_n$ we set the profile rescaled by $1/\sqrt{n}$ as 
\begin{equation}\label{eq:1-2}
[\lambda]^{\sqrt{n}}(x) = \frac{1}{\sqrt{n}} \lambda (\sqrt{n} x), \qquad x\in\mathbb{R}.
\end{equation}
The limit shape of Young diagrams with respect to the Plancherel measure can be described 
in a form of weak law of large numbers as follows. 
An element of 
%
%\begin{equation}\label{eq:1-10}
\begin{align*}
\mathbb{D} = \bigl\{ \omega : \mathbb{R} \longrightarrow\mathbb{R} \,\big|\, 
&|\omega(x) -\omega(y)|\leqq |x-y|, \\ 
&\ \omega(x) =|x| \ (x\leqq a, \;b\leqq x) \text{ for some } 
a\leqq 0 \text{ and } b\geqq 0\bigr\}
\end{align*}
%\end{equation}
%
is called a (centered) continuous diagram. 
The smallest closed interval $[a, b]$ satisfying this condition (possibly a singleton $\{0\}$) 
is denoted by $\mathrm{supp}\,\omega$ (support of $\omega\in\mathbb{D}$). 
Let $\varOmega$ denote the continuous diagram (indeed a $C^1$ curve) with 
$\mathrm{supp}\,\varOmega = [-2, 2]$:  
\begin{equation}\label{eq:1-3}
\varOmega (x) = \begin{cases} \frac{2}{\pi} \bigl( x \arcsin \frac{x}{2} + \sqrt{4-x^2}\bigr), & |x|\leqq 2, \\ 
|x|, & |x|>2. \end{cases}
\end{equation}
Then, $[\lambda]^{\sqrt{n}}$ converges to $\varOmega$ in $\mathbb{D}$ in probability 
$\mathrm{M}_{\mathrm{Pl}}^{(n)}$ as $n\to\infty$. 
Namely, for any $\epsilon >0$, it holds 
\begin{equation}\label{eq:1-4}
\lim_{n\to\infty} \mathrm{M}_{\mathrm{Pl}}^{(n)} \Bigl( \Bigl\{ \lambda\in\mathbb{Y}_n \,\Big|\, 
\sup_{x\in\mathbb{R}} \bigl| [\lambda]^{\sqrt{n}}(x) - \varOmega(x) \bigr| \geqq\epsilon \Bigr\} \Bigr) = 0.
\end{equation}
This result of the limit shape is a static property for the Plancherel ensemble. 
In \cite{Ho15} we treated a dynamical limit shape, in other words, evolution of the limit shapes along 
macroscopic time. 
We considered a continuous time Markov chain keeping the Plancherel measure invariant, took a diffusive 
scaling limit in time and space, and found limit shape (or macroscopic profile) $\omega_t$ depending 
macroscopic time $t$. 
A pioneering work about time evolution of profiles of Young diagrams is done by \cite{FuSa10}. 

The purpose of the present paper is to introduce a pausing time not necessarily obeying an exponential distribution 
instead of sticking to Markovian property for the microscopic dynamics of continuous time and to observe how 
the pausing time distribution produces an effect on the scale of micro-macro correspondence 
and macroscopic evolution of the limit shape. 
Actually, we will see an anomalous effect given by a pausing time distribution with a heavy tail. 

Let us begin with recalling the restriction-induction (Res-Ind) chain on Young diagrams. 
For a finite group $G$ and its subgroup $H$, composing restriction of an irreducible representation of $G$ 
to $H$ ($\mathrm{Res}^G_H$) and induction from $H$ to $G$ ($\mathrm{Ind}_H^G$), 
and counting the dimensions of irreducible decompositions, we get a transition probability on $\widehat{G}$. 
Namely, for $\lambda\in\widehat{G}$ and $\nu\in\widehat{H}$, we have 
\[
\mathrm{Res}^G_H \lambda  \ \cong\ \bigoplus_{\xi\in\widehat{H}} [c_{\lambda\xi}] \xi, 
\qquad 
\mathrm{Ind}^G_H \nu  \ \cong\ \bigoplus_{\mu\in\widehat{G}} [c_{\mu\nu}] \mu
\] 
with multiplicities $c_{\lambda\nu} = [\mathrm{Res}^G_H \lambda, \nu] = [\mathrm{Ind}^G_H \nu, \lambda]$, 
and the transition probabilities 
\begin{align}
&P^\downarrow_{\lambda\nu} = \frac{c_{\lambda\nu}\dim\nu}{\dim\lambda}, \quad 
P^\uparrow_{\nu\mu} = \frac{c_{\mu\nu}\dim\mu}{[G:H]\dim\nu}, & 
&\lambda, \mu\in \widehat{G}, \ \nu\in\widehat{H}, \notag \\ 
&P_{\lambda\mu} =\sum_{\nu\in\widehat{H}} P^\downarrow_{\lambda\nu}P^\uparrow_{\nu\mu}, & 
&\lambda, \mu\in \widehat{G} \label{eq:1-5}
\end{align}
by taking the dimensions of both sides. 
The Plancherel measure on $\widehat{G}$ defined by 
%
%\begin{equation}\label{eq:1-6}
\[ 
\mathrm{M}_{\mathrm{Pl}}^G (\{\lambda\}) = \frac{(\dim\lambda)^2}{|G|}, 
\qquad \lambda\in \widehat{G}
\] 
%\end{equation}
%
makes $P = (P_{\lambda\mu})$ of \eqref{eq:1-5} symmetric: 
%
%\begin{equation}\label{eq:1-7}
\[ 
\mathrm{M}_{\mathrm{Pl}}^G (\{\lambda\}) P_{\lambda\mu} = 
\mathrm{M}_{\mathrm{Pl}}^G (\{\mu\}) P_{\mu\lambda}, \qquad 
\lambda, \ \mu\in \widehat{G}.
\] 
%\end{equation}
%
Specializing in $G=\mathfrak{S}_n$ (the symmetric group of degree $n$) and 
$H=\mathfrak{S}_{n-1}$, and identifying $\widehat{\mathfrak{S}_n}$ with $\mathbb{Y}_n$, 
we get from \eqref{eq:1-5} transition matrix $P^{(n)} = (P_{\lambda\mu})$ of degree 
$|\mathbb{Y}_n|$ which keeps the Plancherel measure on $\mathbb{Y}_n$ invariant. 
Note that in this case 
\[ 
c_{\lambda\nu} = \begin{cases} 1, & \text{if } \nu\nearrow\lambda, \\ 
0, & \text{otherwise} \end{cases}
\] 
where $\nu\nearrow\lambda$ indicates that $\nu$ is formed by removing a box of $\lambda$. 
The Markov chain determined by $P^{(n)}$ is the Res-Ind chain on $\mathbb{Y}_n$. 
In this chain, a one step transition admits non-local movement of a corner box in a Young diagram. 
The Res-Ind chain was treated in \cite{Fu04}, \cite{Fu05}, and \cite{BoOl09}. 

Let us construct a continuous time random walk on $\mathbb{Y}_n$, not necessarily Markovian, from 
transition matrix $P^{(n)}$. 
We mention \cite{We94} as a nice reference on such a non-Markovian continuous time random walk. 
Consider Markov chain $(Z^{(n)}_k)_{k\in\{0,1,2,\cdots\}}$ on $\mathbb{Y}_n$ having 
transition matrix $P^{(n)}$ and initial distribution $M^{(n)}_0$. 
Let $(\tau_j)_{j\in\mathbb{N}}$ be i.i.d. random variables independent also of $(Z^{(n)}_k)$, 
each obeying $\psi(dx)$ on $[0, \infty)$. 
This sequence yields counting process $(N_s)_{s\geqq 0}$ in which pausing intervals are given 
by $\tau_j$'s: 
\[ 
N_s = \begin{cases} j, & \tau_1+ \cdots + \tau_j < s\leqq \tau_1 + \cdots +\tau_{j+1} \\ 
0, & s\leqq \tau_1 \end{cases}, \qquad 
N_0 =0 \quad \text{a.s.}
\] 
Note that we assume nontriviality of $\psi$, $\psi\bigl( (0, \infty)\bigr) >0$, which implies 
$\tau_1+\cdots+\tau_j$ diverges to $\infty$ a.s. as $j\to\infty$. 
Set 
\begin{equation}\label{eq:1-71}
X^{(n)}_s = Z^{(n)}_{N_s}, \qquad s\geqq 0. 
\end{equation}
The process $(X^{(n)}_s)_{s\geqq 0}$ is a desired continuous time random walk on $\mathbb{Y}_n$. 
We have 
\begin{align*}
&\mathrm{Prob}(X^{(n)}_s =\mu \,|\, X^{(n)}_0 =\lambda) \notag \\ 
&= 
\sum_{j=0}^\infty \mathrm{Prob}(Z^{(n)}_{N_s} =\mu \,|\, N_s =j, Z^{(n)}_0=\lambda) \, 
\mathrm{Prob}(N_s=j \,|\, Z^{(n)}_0 =\lambda) \notag \\ 
&= 
\sum_{j=0}^\infty \mathrm{Prob}(Z^{(n)}_j =\mu \,|\, Z^{(n)}_0=\lambda) \, 
\mathrm{Prob}(\tau_1+\cdots +\tau_j\leqq s, \tau_1+\cdots +\tau_{j+1} >s) \notag \\ 
&= 
\sum_{j=0}^\infty (P^{(n)\,j})_{\lambda\mu} \int_{[0, s]} \psi\bigl( (s-u, \infty)\bigr) 
\psi^{\ast j}(du)
%\label{eq:1-8}
\end{align*}
where $\psi^{\ast j}$ means the ordinary $j$-fold convolution power of $\psi$. 
Regarding initial distribution 
$M^{(n)}_0(\{\lambda\}) = \mathrm{Prob}(X^{(n)}_0 =\lambda)$ 
as a row vector of degree $|\mathbb{Y}_n|$, we have the distribution at time $s$ as 
\begin{equation}\label{eq:1-81}
M^{(n)}_s (\{\mu\}) = \mathrm{Prob}(X^{(n)}_s =\mu) = 
\sum_{j=0}^\infty (M^{(n)}_0 P^{(n)\, j})_\mu 
\int_{[0, s]} \psi\bigl( (s-u, \infty)\bigr) \psi^{\ast j}(du).
\end{equation}

In \eqref{eq:1-4} we stated a result of the limit shape for a sequence of the Plancherel measures 
$\{(\mathbb{Y}_n, \mathrm{M}_{\mathrm{Pl}}^{(n)})\}_{n\in\mathbb{N}}$. 
In more general ensembles, we intend to observe concentration of rescaled Young diagrams 
$[\lambda]^{\sqrt{n}}$ ($\lambda\in\mathbb{Y}_n$) at a continuous diagram $\omega$ as $n\to\infty$. 
Taking into account the algebraic structure for functions of the coordinates of Young diagrams, 
we recognize that a formulation under a stronger convergence than the weak law of large numbers 
with respect to the uniform topology on $\mathbb{D}$ 
as in \eqref{eq:1-4} is more suitable for our purpose, as follows. 
Set 
$\mathbb{D}^{(a)} = \{ \omega\in\mathbb{D} \,|\, \mathrm{supp}\,\omega \subset [-a, a]\}$ 
for $a>0$ to have 
$\mathbb{D} = \bigcup_{a>0} \mathbb{D}^{(a)}$. 
By using the $k$th moment $M_k(\mathfrak{m}_\omega)$ of transition measure $\mathfrak{m}_\omega$ 
of $\omega\in\mathbb{D}$, let us equip $\mathbb{D}$ with the topology induced by the family 
of pseudo-distances 
$\{ | M_k(\mathfrak{m}_{\omega_1}) - M_k(\mathfrak{m}_{\omega_2}) | \}_{k\in\mathbb{N}}$ 
and call it the moment topology on $\mathbb{D}$. 
The moment topology and the uniform one are equivalent on $\mathbb{D}^{(a)}$. 
We equip $\mathbb{D}^{(a)}$ with this topology. 
Then $\mathbb{D}$ has the topologies in a stronger order: inductive limit topology of 
$\{\mathbb{D}^{(a)}\}$, moment topology, and uniform topology.  
See \cite[Sections~3.1 and 3.3]{Ho16}. 
We thus have a stronger condition if the convergence of \eqref{eq:1-4} is replaced by 
\[ 
\lim_{n\to\infty} \mathbb{E}_{\mathrm{M}_{\mathrm{Pl}}^{(n)}} \bigl[ 
\bigl( M_k(\mathfrak{m}_{[\lambda]^{\sqrt{n}}}) - M_k(\mathfrak{m}_\varOmega)\bigr)^2 
\bigr] = 0, \qquad k\in\mathbb{N}
\] 
where $\mathbb{E}_M$ denotes the expectation in variable $\lambda$ under probability $M$. 
Furthermore, to make algebraic arguments transparent, we consider the convergence of all mixed 
moments without restricting such second order ones. 
This brings us to the following formulation.

\begin{definition}\label{def:1-0}
Let $\{ (\mathbb{Y}_n, M^{(n)})\}_{n\in\mathbb{N}}$ be a sequence of probability spaces. 
If there exists $\omega\in\mathbb{D}$ such that 
\begin{equation}\label{eq:1-8-5}
\lim_{n\to\infty} \mathbb{E}_{M^{(n)}} \bigl[ M_{k_1}(\mathfrak{m}_{[\lambda]^{\sqrt{n}}}) 
\cdots M_{k_p}(\mathfrak{m}_{[\lambda]^{\sqrt{n}}}) \bigr] = 
M_{k_1}(\mathfrak{m}_\omega) \cdots M_{k_p}(\mathfrak{m}_\omega)
\end{equation}
holds for any $p\in\mathbb{N}$ and any $k_1, \cdots, k_p \in \{2,3,\cdots\}$, we say that 
$\{ (\mathbb{Y}_n, M^{(n)})\}_{n\in\mathbb{N}}$ admits concentration at $\omega\in\mathbb{D}$ 
as $n\to\infty$.
\end{definition}

Definition~\ref{def:1-0} obviously yields the weak law of large numbers 
with respect to the uniform topology on $\mathbb{D}$. 
A mechanism causing such a concentration phenomenon for a group-theoretical ensemble 
$\{ (\mathbb{Y}_n, M^{(n)})\}_{n\in\mathbb{N}}$ was pointed out by Biane \cite{Bi01} 
as approximate factorization property. 
Approximate factorization property can be described in several equivalent ways. 
Here we define it in terms of irreducible characters of the symmetric groups as follows. 
The irreducible character of $\mathfrak{S}_n$ corresponding to $\lambda\in\mathbb{Y}_n$ 
is denoted by $\chi^\lambda$. 
The value it takes at an element of the conjugacy class of $\mathfrak{S}_n$ corresponding to 
$\rho\in\mathbb{Y}_n$ is $\chi^\lambda_\rho$. 
Normalization of $\chi^\lambda$ yields $\widetilde{\chi}^\lambda = \chi^\lambda /\dim\lambda$. 
Set $\mathbb{Y}^\times = \{ \lambda\in\mathbb{Y} \,|\, m_1(\lambda) =0\}$. 
When we fix a type of a conjugacy class and let the size $n$ tend to infinity, we use a 
convenient notation as 
\[ 
(\rho, 1^{n-|\rho|}) = \rho \sqcup (1^{n-|\rho|}), \qquad \rho\in\mathbb{Y}^\times
\] 
for the Young diagram of size $n$ indicating a type of a conjugacy class. 

\begin{definition}\label{def1-1}
A sequence of probability spaces $\{ (\mathbb{Y}_n, M^{(n)})\}_{n\in\mathbb{N}}$ 
is said to satisfy approximate factorization property if 
\begin{equation}\label{eq:1-9}
\mathbb{E}_{M^{(n)}}[ 
\widetilde{\chi}^\lambda_{(\rho\sqcup\sigma, 1^{n-|\rho|-|\sigma|})} ] - 
\mathbb{E}_{M^{(n)}}[  
\widetilde{\chi}^\lambda_{(\rho, 1^{n-|\rho|})} ] \, 
\mathbb{E}_{M^{(n)}}[ 
\widetilde{\chi}^\lambda_{(\sigma, 1^{n-|\sigma|})} ] 
= o \bigl( n^{-\frac{1}{2}(|\rho|-l(\rho)+|\sigma|-l(\sigma))} \bigr) 
\end{equation}
as $n\to\infty$ holds for any $\rho, \sigma\in \mathbb{Y}^\times$.
\end{definition}
Concerning the decay order in the right hand side of \eqref{eq:1-9}, see also \eqref{eq:2-2} in Section~2. 
Expectations of irreducible characters seen in \eqref{eq:1-9} are analogous objects to 
characteristic functions of probabilities. 
Since \eqref{eq:1-9} says that characteristic functions are nearly factorizable along cycle decomposition 
with small error terms in some sense, 
approximate factorization property is regarded as an analogous, but much weaker, notion to independence. 
Applying approximate factorization property, Biane extended the concentration phenomenon 
\eqref{eq:1-4} for the Plancherel measure to a wide variety of interesting models in \cite{Bi01}. 
For convenience of later reference, we here give a statement in the following form. 
See also Section~4.4 of \cite{Ho16} for a proof in detail. 

\begin{proposition}\label{prop:1-1}
A sequence of probability spaces $\{(\mathbb{Y}_n, M^{(n)})\}_{n\in\mathbb{N}}$ 
admits concentration at $\omega\in\mathbb{D}$ (in the sense of Definition~\ref{def:1-0}) 
if and only if 

{\rm (i)} it satisfies approximate factorization property \eqref{eq:1-9}, 

{\rm (ii)} the limit of the expected value at $j$-cycle $(j, 1^{n-j})$ 
\begin{equation}\label{eq:1-11}
\lim_{n\to\infty} n^{\frac{j-1}{2}} \mathbb{E}_{M^{(n)}} [  
\widetilde{\chi}^{\,\cdot}_{(j, 1^{n-j})}] = r_{j+1}, 
\qquad j\in \{2,3,\cdots\}
\end{equation}
exists and has an order of at most $j$th power: 
\begin{equation}\label{eq:1-111}
|r_j| \leqq b^j \qquad \text{for some } b>0.
\end{equation}
%
%Then we have concentration at a continuous diagram $\omega\in\mathbb{D}$; 
%namely, $[\lambda]^{\sqrt{n}}$ converges to $\omega$ in $\mathbb{D}$ in probability $M^{(n)}$ 
%as $n\to\infty$. 
%for any $\epsilon >0$ 
%%
%\begin{equation}\label{eq:1-12}
%\lim_{n\to\infty} M^{(n)} \Bigl( \Bigl\{ \lambda\in\mathbb{Y}_n \,\Big|\, 
%\sup_{x\in\mathbb{R}} \bigl| [\lambda]^{\sqrt{n}}(x) - \omega(x) \bigr| \geqq\epsilon \Bigr\} \Bigr) = 0.
%\end{equation}
%%
In this situation, the limit shape $\omega$ is characterized by free cumulants of its transition measure 
$\mathfrak{m}_\omega$ as 
%
%\begin{equation}\label{eq:1-13}
\[ 
R_1(\mathfrak{m}_\omega) =0, \quad R_2(\mathfrak{m}_\omega) =1, \quad 
R_{j+1}(\mathfrak{m}_\omega) = r_{j+1} \ (j\in\{2,3,\cdots\}).
\] 
%\end{equation}
%
\end{proposition}
A procedure of computing $\omega$ from a sequence of free cumulants 
$\{R_j(\mathfrak{m}_\omega)\}_{j\in\mathbb{N}}$ is given by the Markov transform 
(see Appendix). 

%In what follows, we assume, though somewhat technically, that pausing time distribution $\psi$ is 
%absolutely continuous with $\psi(x)$ as the density. 
Let $\varphi$ be the characteristic function (Fourier transform) of $\psi$: 
%
%\begin{equation}\label{eq:1-14}
\[ 
\varphi(\xi) = \int_{[0, \infty)} e^{i\xi x}\psi(dx),  
%= \int_0^\infty e^{i\xi x}\psi(x)dx, 
\qquad \xi\in\mathbb{R}.
\] 
%\end{equation}
%
%We see $\lim_{\xi\to\pm\infty} \varphi(\xi) =0$ from the Riemann--Lebesgue theorem. 
Differentiability at $\xi =0$ of $\varphi$ follows if $\psi$ has the mean. 

The first result of the present paper is the following scaling limit of the continuous time 
random walk $(X^{(n)}_s)_{s\geqq 0}$. 

\begin{theorem}\label{th:1-1} 
Let $(X^{(n)}_s)_{s\geqq 0}$ be the continuous time random walk of \eqref{eq:1-71}. 
For any microscopic time $s\geqq 0$, let the distribution at time $s$ be 
\[ 
M^{(n)}_s(\{\lambda\}) = \mathrm{Prob}(X^{(n)}_s =\lambda), \qquad \lambda\in\mathbb{Y}_n. 
\] 
Assume that the sequence of initial distributions 
$\{(\mathbb{Y}_n, M^{(n)}_0)\}_{n\in\mathbb{N}}$ admits concentration at 
$\omega_0 \in \mathbb{D}$ in the sense of Definition~\ref{def:1-0}. 
Assume also that the pausing time distribution $\psi$ has the mean $m$ and that 
the characteristic function $\varphi$ of $\psi$ satisfies the integrability condition 
\begin{equation}\label{eq:1-15}
\int_{\{|\xi|\geqq \delta\}} \Bigl| \frac{\varphi(\xi)}{\xi}\Bigr| d\xi < \infty 
\qquad \text{for some } \delta >0.
\end{equation}
Then, by considering $s=tn$ for macroscopic time $t>0$, 
$\{(\mathbb{Y}_n, M^{(n)}_{tn})\}_{n\in\mathbb{N}}$ inherits 
%approximate factorization property together with \eqref{eq:1-11} and \eqref{eq:1-111}, 
the condition of Definition~\ref{def:1-0} and hence admits 
concentration at some $\omega_t \in\mathbb{D}$. 
%namely $[\lambda]^{\sqrt{n}}$ converges to $\omega_t$ in $\mathbb{D}$ in probability 
%$M^{(n)}_{tn}$ as $n\to\infty$. 
%for any $\epsilon >0$, 
%%
%\begin{equation}\label{eq:1-16}
%\lim_{n\to\infty} M^{(n)}_t \Bigl( \Bigl\{ \lambda\in\mathbb{Y}_n \,\Big|\, 
%\sup_{x\in\mathbb{R}} \bigl| [\lambda]^{\sqrt{n}}(x) - \omega_t(x)\bigr| \geqq\epsilon 
%\Bigr\} \Bigr) =0.
%\end{equation}
%%
The transition measure of the limit shape $\omega_t$ is given by 
\begin{equation}\label{eq:1-17}
\mathfrak{m}_{\omega_t} = (\mathfrak{m}_{\omega_0})_{e^{-t/m}} \boxplus 
(\mathfrak{m}_\varOmega)_{1-e^{-t/m}}, \qquad t>0.
\end{equation}
In \eqref{eq:1-17}, $(\;\cdot\;)_c$ denotes free compression of rank $c$, $\boxplus$ denotes 
free convolution, and $\varOmega$ is the limit shape \eqref{eq:1-3} of Vershik--Kerov and Logan--Shepp. 
Equivalently to \eqref{eq:1-17} in terms of the free cumulants, we have 
\begin{equation}\label{eq:1-18}
R_1(\mathfrak{m}_{\omega_t})=0, \quad R_2(\mathfrak{m}_{\omega_t})=1, \qquad 
R_{k+1}(\mathfrak{m}_{\omega_t}) = R_{k+1}(\mathfrak{m}_{\omega_0})e^{-kt/m} 
\quad (k\geqq 2).
\end{equation}
\end{theorem}

We note 
it is possible to choose a desired sequence of initial distributions for arbitrarily prescribed $\omega_0\in\mathbb{D}$ 
such that $\int_{\mathbb{R}} (\omega_0(x) -|x|)dx =2$. 
We see from \eqref{eq:1-17} 
\[ 
\int_{\mathbb{R}} \bigl( \omega_t(x)-|x|\bigr) dx =2 \qquad \text{and} \qquad  
\lim_{t\to\infty} \omega_t = \varOmega \  \text{in} \  \mathbb{D}.
\] 

A main result in \cite{Ho15} is a special case of Theorem~\ref{th:1-1}, in which $(X^{(n)}_s)$ 
is a continuous time Markov chain, or the pausing time obeys an exponential distribution (with 
mean $1$). 
Properties of such free convolution with semi-circular distributions as \eqref{eq:1-17} were 
treated in detail in \cite{Bi97}. 
See \cite{NiSp06} and Appendix also for necessary notions in free probability theory. 
Proof of Theorem~\ref{th:1-1} is presented in Section~2. 
In the situation of Theorem~\ref{th:1-1}, microscopic time $s=tn$ is of order $n$ while the rescale 
of space is of $\sqrt{n}$ as in \eqref{eq:1-2}. We thus took a diffusive scaling limit. 
The Stieltjes transform of $\mathfrak{m}_{\omega_t}$ 
\[ 
G(t, z) = \int_{\mathbb{R}} \frac{1}{z-x} \mathfrak{m}_{\omega_t}(dx), 
\qquad z\in\mathbb{C}^+
\] 
satisfies the partial differential equation 
\begin{equation}\label{eq:1-18.5}
m \frac{\partial G}{\partial t}(t, z) = - G(t, z) \frac{\partial G}{\partial z}(t, z) + 
\frac{1}{G(t, z)} \frac{\partial G}{\partial z}(t, z) + G(t, z), 
\end{equation}
which is derived from \cite[Theorem~3.3]{Ho15}.

It is remarkable that Biane \cite{Bi98} observed the appearance of free convolution 
and free compression as a result of concentration of rescaled Young diagrams in static models 
which are produced by irreducible decomposition of induction (outer product) and restriction 
for irreducible representations of symmetric groups. 
Although we do not use these results in the present paper, the structure of \eqref{eq:1-17} 
at each macroscopic time $t>0$ seems to be natural in this framework.

On the other hand, when we consider the case where a microscopic pausing time distribution is 
heavy-tailed so as not to have the mean any more, it is naturally expected that limiting behavior will be 
different from the one in Theorem~\ref{th:1-1}. 
The second result of the present paper illustrates such an observation. 
Let us take a pausing time obeying the one-sided stable distribution $\psi$ of exponent $\alpha\in (0, 1)$ 
whose characteristic function is given by 
\begin{equation}\label{eq:1-20}
\varphi(\xi) = e^{-|\xi|^\alpha (1- i \tan(\pi\alpha /2) \mathrm{sgn}(\xi))}, \qquad \xi\in\mathbb{R}.
\end{equation}
The distribution $\psi$ is absolutely continuous. 
Especially in the simplest case of the exponent $1/2$, its density is expressed as 
%
%\begin{equation}\label{eq:1-19}
\[ 
\frac{1}{\sqrt{2\pi}} x^{-\frac{3}{2}} e^{-\frac{1}{2x}} 1_{(0, \infty)}(x).
\] 
%\end{equation}
%
See e.g. \cite{Lu60} for one-sided stable distributions and their characteristic functions. 
As for the scaling limit for continuous time random walk $(X^{(n)}_s)_{s\geqq 0}$ with such a $\psi$ 
as its pausing time distribution, it proves that approximate factorization property of an initial ensemble 
is not propagated along positive macroscopic time.   

\begin{theorem}\label{th:1-2} 
Let $(X^{(n)}_s)_{s\geqq 0}$ be the continuous time random walk of \eqref{eq:1-71}. 
For any microscopic time $s\geqq 0$, let the distribution at time $s$ be 
\[ 
M^{(n)}_s(\{\lambda\}) = \mathrm{Prob}(X^{(n)}_s =\lambda), \qquad \lambda\in\mathbb{Y}_n. 
\] 
Assume that the sequence of initial distributions 
$\{(\mathbb{Y}_n, M^{(n)}_0)\}_{n\in\mathbb{N}}$ satisfies approximate factorization property 
together with \eqref{eq:1-11} and 
\eqref{eq:1-111}, and hence has concentration at $\omega_0 \in \mathbb{D}$. 
Assume also \eqref{eq:1-20} for the pausing time distribution. 
For macroscopic time $t>0$, let $s = t\theta_n$, where the scaling factor 
$\theta : \mathbb{N} \longrightarrow \mathbb{R}_+$ is taken to be such that 
\[ 
\text{\rm (i)} \ \frac{\theta_n}{n^{1/\alpha}}\ \longrightarrow \ 0, \quad 
\text{\rm (ii)} \ \frac{\theta_n}{n^{1/\alpha}}\ \longrightarrow \ \infty, \quad 
\text{\rm (iii)} \ \frac{\theta_n}{n^{1/\alpha}}\ \longrightarrow \ 1 \quad \text{as} \quad n\to\infty.
\] 
Then, in either case of {\rm (i)} or {\rm (ii)}, 
$\{(\mathbb{Y}_n, M^{(n)}_{t\theta_n})\}_{n\in\mathbb{N}}$ inherits 
approximate factorization property together with \eqref{eq:1-11} and \eqref{eq:1-111}, and hence has 
concentration at $\omega_t \in\mathbb{D}$. 
The limit shape $\omega_t$ is, however, rather trivial so that 

{\rm (i)} \ $\omega_t = \omega_0$, that is, no macroscopic evolution observed 

{\rm (ii)}\ $\omega_t = \varOmega$ for any $t>0$, that is, macroscopic evolution completed at once. 

\noindent 
In the case of {\rm (iii)}, we have the convergence of the averaged quantities 
\begin{multline}\label{eq:1-21}
\lim_{n\to\infty} \mathbb{E}_{M^{(n)}_{t\theta_n}} 
\bigl[ R_{k+1}(\mathfrak{m}_{[\lambda]^{\sqrt{n}}}) \bigr] \\ 
= R_{k+1}(\mathfrak{m}_{\omega_0}) \, 
\frac{\sin \pi\alpha}{\pi\alpha} \int_0^\infty e^{-t(k\xi\cos(\pi\alpha/2))^{1/\alpha}}\, 
\frac{d\xi}{\xi^2+ 2\xi\cos(\pi\alpha)+1}, \qquad k\geqq 2 
\end{multline}
for any initial $\omega_0\in\mathbb{D}$. 
However, $\{(\mathbb{Y}_n, M^{(n)}_{t\theta_n})\}_{n\in\mathbb{N}}$ inherits 
approximate factorization property if and only if $\omega_0 = \varOmega$ 
(hence $\omega_t = \varOmega$ for any $t$ also). 
%
%In terms of free cumulants, 
%%
%\begin{equation}\label{eq:1-21}
%R_1(\mathfrak{m}_{\omega_t})=0, \quad R_2(\mathfrak{m}_{\omega_t})=1, \qquad 
%R_{k+1}(\mathfrak{m}_{\omega_t}) = R_{k+1}(\mathfrak{m}_{\omega_0}) c_k(t) 
%\quad (k\geqq 2)
%\end{equation}
%
%
%\begin{equation}\label{eq:1-22}
%\gamma (u) = \frac{2}{\pi} \int_0^\infty e^{-\frac{u \xi^2}{2}} \frac{1}{\xi^2 +1} d\xi,
%\qquad u\geqq 0.
%\end{equation}
%
\end{theorem}
Proof of Theorem~\ref{th:1-2} is given in Section~2. 
%Asymptotically $c_k(t)$ in \eqref{eq:1-22} ($k\in\mathbb{N}$) satisfies 
%%
%\begin{equation}\label{eq:1-23}
%c_k(t) \ \sim \ \sqrt{\frac{2}{\pi}} \frac{1}{k\sqrt{t}} \quad \text{as} \quad t\to\infty.
%\end{equation}
%%
%Hence the behavior along macroscopic time $t$ is different from \eqref{eq:1-18} in Theorem~\ref{th:1-1}.

The subsequent sections are organized as follows. 
In Section~2, we give proofs of the theorems. 
However, proofs of the essential propositions involving computational details are postponed until Section~3. 
Our method relies on Fourier analysis (both classical and more group-theoretical). 
Usefulness of Fourier analysis is already suggested in \cite{We94} in treating continuous time random walks 
under general pausing time. 
%Remark contains some comments. 
%Appendix collects brief reviews of some necessary notions.

%\subsection{A subsection}

\section{Proofs of Theorems}

The mechanism of propagating approximate factorization property along macroscopic time is exactly the same 
as treated in \cite{Ho15}. 
See also Section~5.2 in \cite{Ho16} for more information. 

Normalizing an irreducible character of a symmetric group, let us consider a function on $\mathbb{Y}$ 
for each $\rho\in\mathbb{Y}$
%
%\begin{equation}\label{eq:2-1}
\[ 
\Sigma_\rho(\lambda) = \begin{cases} |\lambda|^{\downarrow |\rho|}\, 
\widetilde{\chi}^\lambda_{(\rho, 1^{|\lambda|-|\rho|})}, & |\lambda|\geqq |\rho|, \\ 
0, & |\lambda|< |\rho|, \end{cases} \qquad\quad \lambda\in\mathbb{Y}
\] 
%\end{equation}
%
where $n^{\downarrow k} = n(n-1)\cdots (n-k+1)$. 
The algebra $\mathbb{A}$ consisting of all linear hulls of $\Sigma_\rho$'s plays a fundamental role 
in the dual approach due to Kerov and Olshanski. 
Basically, our harmonic analysis is developed in this algebra. 
See \cite{IvOl02} for its structure. 
For a sequence of probability spaces $\{(\mathbb{Y}_n, M^{(n)})\}_{n\in\mathbb{N}}$, 
approximate factorization property \eqref{eq:1-9} and \eqref{eq:1-11} are rephrased in terms of 
$\Sigma_\rho$'s:  
\begin{equation}\label{eq:2-2}
\mathbb{E}_{M^{(n)}}[ \Sigma_{\rho\sqcup\sigma} ] - 
\mathbb{E}_{M^{(n)}}[ \Sigma_\rho ] \, \mathbb{E}_{M^{(n)}}[ \Sigma_\sigma ] = 
o\bigl( n^{\frac{1}{2}(|\rho|+l(\rho)+|\sigma|+l(\sigma))}\bigr) 
\end{equation}
as $n\to\infty$ for $\rho, \sigma\in\mathbb{Y}^\times$, and 
\begin{equation}\label{eq:2-3}
\lim_{n\to\infty} n^{-\frac{j+1}{2}} \mathbb{E}_{M^{(n)}}[ \Sigma_j ] = r_{j+1} 
\end{equation}
for $j\in\{2,3,\cdots\}$. 
We note also that \eqref{eq:2-2} and \eqref{eq:2-3} yield 
\begin{equation}\label{eq:2-31}
\mathbb{E}_{M^{(n)}}[ \Sigma_\rho] = O \bigl( n^{\frac{1}{2}(|\rho|+l(\rho))}\bigr), 
\qquad \rho\in\mathbb{Y}^\times.
\end{equation}
Let $\mathrm{wt}$ denote the weight degree in $\mathbb{A}$. 
Since $\mathrm{wt}(\Sigma_j) = j+1$ holds, the right hand side of \eqref{eq:2-2} is 
$o(\sqrt{n}^{\,\mathrm{wt}(\Sigma_\rho)+\mathrm{wt}(\Sigma_\sigma)})$. 
We know the relation in $\mathbb{A}$ 
\begin{equation}\label{eq:2-4}
\Sigma_k(\lambda) = R_{k+1}(\mathfrak{m}_\lambda) + 
\ll \text{lower terms with } \mathrm{wt} \leqq k-1 \gg, \qquad k\in\mathbb{N}, 
\end{equation}
which is a decisive formula connetcing irreducible characters of symmetric groups with transition measures 
of Young diagrams. 
Actually, \eqref{eq:2-4} makes our scaling arguments tranparent. 
The right hand side of \eqref{eq:2-4} is a polynomial, known as a Kerov polynomial, in 
$R_j(\mathfrak{m}_\lambda)$'s.

The following formula for transition matrix $P^{(n)}$ of the Res-Ind chain is a key observation about 
propagation of approximate factorization property. 
Regarding $\Sigma_\rho|_{\mathbb{Y}_n}$ as the column vector consisting of the values of 
$\Sigma_\rho$ on $\mathbb{Y}_n$, we have 
\begin{equation}\label{eq:2-5}
P^{(n)} \,\Sigma_\rho|_{\mathbb{Y}_n} = \Bigl( 1- \frac{|\rho|-m_1(\rho)}{n} \Bigr) 
\Sigma_\rho|_{\mathbb{Y}_n}, \qquad \rho\in\mathbb{Y}.
\end{equation}
Formula \eqref{eq:2-5} is obtained by the induced character formula. 

Let $M^{(n)}_s$ denote the distribution of continuous time random walk $(X^{(n)}_s)_{s\geqq 0}$ 
at time $s$. 
For $\rho\in\mathbb{Y}$, \eqref{eq:1-81} and \eqref{eq:2-5} yield 
\begin{align}
\mathbb{E}_{M^{(n)}_s}[ \Sigma_\rho] &= 
\sum_{\mu\in\mathbb{Y}_n} \Bigl( \sum_{j=0}^\infty (M^{(n)}_0 P^{(n)\,j})_\mu 
\int_{[0, s]} \psi \bigl( (s-u, \infty)\bigr) \psi^{\ast j}(du) \Bigr) \Sigma_\rho (\mu) \notag \\  
&= 
\sum_{j=0}^\infty M^{(n)}_0 P^{(n)\,j} \Sigma_\rho|_{\mathbb{Y}_n} 
\int_{[0, s]} \psi \bigl( (s-u, \infty)\bigr) \psi^{\ast j}(du) \notag \\ 
&= 
\Bigl( \sum_{j=0}^\infty \bigl( 1- \frac{|\rho|-m_1(\rho)}{n}\bigr)^j 
\int_{[0, s]} \psi \bigl( (s-u, \infty)\bigr) \psi^{\ast j}(du) \Bigr) 
\mathbb{E}_{M^{(n)}_0} [\Sigma_\rho].
\label{eq:2-6}
\end{align}
Especially in \eqref{eq:2-6}, considering a $k$-cycle, set 
\begin{equation}\label{eq:2-7}
f (k,n,s) = \sum_{j=0}^\infty \bigl( 1-\frac{k}{n}\bigr)^j \int_{[0, s]} \psi \bigl( (s-u, \infty)\bigr) \psi^{\ast\,j}(du) , 
\qquad k\geqq 2.
\end{equation}

\begin{proposition}\label{prop:2-1}
Under the assumptions and notations for the pausing time distribution 
in Theorem~\ref{th:1-1}, let $s=tn$ in \eqref{eq:2-7}. 
Then we have 
\begin{equation}\label{eq:2-8}
\lim_{n\to\infty} f(k, n, tn) = e^{-kt/m}, \qquad k\in\{2,3,\cdots\}.
\end{equation}
\end{proposition}

\begin{proposition}\label{prop:2-2}
Under the assumptions and notations for the pausing time distribution 
in Theorem~\ref{th:1-2}, we have for $k\in\{2,3,\cdots\}$ 
\begin{align}
&\lim_{n\to\infty} f(k,n, t\theta_n) = 1 & &\text{if } \lim_{n\to\infty} \theta_n /n^{1/\alpha} =0, 
\label{eq:2-8} \\ 
&\lim_{n\to\infty} f(k,n, t\theta_n) = 0 & &\text{if } \lim_{n\to\infty} \theta_n /n^{1/\alpha} =\infty, 
\label{eq:2-9} \\ 
&\lim_{n\to\infty} f(k,n, t\theta_n) = \frac{\sin \pi\alpha}{\pi\alpha} \int_0^\infty 
\frac{e^{-t(k\xi\cos(\pi\alpha/2))^{1/\alpha}}}{\xi^2+ 2\xi\cos(\pi\alpha)+1}\, d\xi & 
&\text{if } \lim_{n\to\infty} \theta_n /n^{1/\alpha} =1 
\label{eq:2-10}
\end{align}
according to the cases of {\rm (i)}, {\rm (ii)} and {\rm (iii)} in Theorem~\ref{th:1-2}.
\end{proposition}
Proofs of Propositions~\ref{prop:2-1} and \ref{prop:2-2} are given in Section~3. 
Let us complete the proofs of Theorems~\ref{th:1-1} and \ref{th:1-2} by using 
Propositions~\ref{prop:2-1} and \ref{prop:2-2}.

\subsection{Proof of Theorem~\ref{th:1-1}}

Let us verify the sequence $\{(\mathbb{Y}_n, M^{(n)}_{tn})\}_{n\in\mathbb{N}}$ satisfies 
\eqref{eq:2-2} ($\Leftrightarrow$ \eqref{eq:1-9}). 
For $\rho,\sigma \in\mathbb{Y}^\times$, \eqref{eq:2-6} and \eqref{eq:2-7} yield
\begin{align}
&\mathbb{E}_{M^{(n)}_{tn}}[ \Sigma_{\rho\sqcup\sigma}] - 
\mathbb{E}_{M^{(n)}_{tn}}[ \Sigma_\rho] \, 
\mathbb{E}_{M^{(n)}_{tn}}[ \Sigma_\sigma] \notag \\ 
&= f(|\rho|+|\sigma|, n, tn) \mathbb{E}_{M^{(n)}_0}[ \Sigma_{\rho\sqcup\sigma}] - 
f(|\rho|, n, tn) f(|\sigma|, n, tn) 
\mathbb{E}_{M^{(n)}_0}[ \Sigma_\rho] \, \mathbb{E}_{M^{(n)}_0}[ \Sigma_\sigma] \notag \\ 
&= \bigl( f(|\rho|+|\sigma|, n, tn) - e^{-(|\rho|+|\sigma|)t/m} \bigr) \, 
\mathbb{E}_{M^{(n)}_0}[ \Sigma_{\rho\sqcup\sigma}] \notag \\ 
&\quad + \bigl( e^{-|\rho|t/m} - f(|\rho|, n, tn)\bigr) \, f(|\sigma|, n, tn) \, 
\mathbb{E}_{M^{(n)}_0}[ \Sigma_\rho] \, \mathbb{E}_{M^{(n)}_0}[ \Sigma_\sigma] \notag \\ 
&\quad + e^{-|\rho|t/m} \bigl( e^{-|\sigma|t/m} - f(|\sigma|, n, tn)\bigr) \, 
\mathbb{E}_{M^{(n)}_0}[ \Sigma_\rho] \, \mathbb{E}_{M^{(n)}_0}[ \Sigma_\sigma] \notag \\ 
&\quad + e^{-(|\rho|+|\sigma|)t/m} \mathbb{E}_{M^{(n)}_0}[ \Sigma_{\rho\sqcup\sigma}] - 
e^{-|\rho|t/m} e^{-|\sigma|t/m} \mathbb{E}_{M^{(n)}_0}[ \Sigma_\rho] \, 
\mathbb{E}_{M^{(n)}_0}[ \Sigma_\sigma] \notag \\ 
&= o \bigl(n^{\frac{1}{2}(|\rho|+l(\rho)+|\sigma|+l(\sigma))}\bigr)
\label{eq:2-1-1}
\end{align}
by taking into account Proposition~\ref{prop:2-1} 
with \eqref{eq:2-2} and \eqref{eq:2-31} for $M^{(n)}_0$. 

To verify \eqref{eq:2-3} ($\Leftrightarrow$ \eqref{eq:1-11}) for $\{(\mathbb{Y}_n, M^{(n)}_{tn})\}_{n\in\mathbb{N}}$, 
we see from \eqref{eq:2-6} and \eqref{eq:2-7} 
\begin{equation}\label{eq:2-1-2}
n^{-\frac{j+1}{2}} \mathbb{E}_{M^{(n)}_{tn}} [\Sigma_j] = f(j, n, tn) \, 
n^{-\frac{j+1}{2}} \mathbb{E}_{M^{(n)}_0} [\Sigma_j] \ \xrightarrow[\;n\to\infty\;] \ 
e^{-jt/m} r_{j+1}
\end{equation}
for $j\in\{2,3,\cdots\}$ by \eqref{eq:1-11} for $M^{(n)}_0$ and Proposition~\ref{prop:2-1}.

We see $M^{(n)}_{tn}$ also satisfies \eqref{eq:1-111} since 
\[ 
|e^{-(j-1)t/m} r_j | \leqq e^{t/m} (e^{-t/m}b)^j 
\] 
holds under \eqref{eq:1-111} for $M^{(n)}_0$.

Finally we look at the free cumulants of transition measure $\mathfrak{m}_{\omega_t}$ of 
\[ 
\omega_t = \lim_{n\to\infty} [\lambda]^{\sqrt{n}} \qquad (\text{ in probability } M^{(n)}_{tn}).
\] 
The first and second ones hold before taking limit. 
We see from \eqref{eq:2-4} and \eqref{eq:2-1-2} 
\begin{equation}\label{eq:2-1-3}
R_{k+1}(\mathfrak{m}_{\omega_t}) = 
\lim_{n\to\infty} n^{-\frac{k+1}{2}} \mathbb{E}_{M^{(n)}_{tn}}[ R_{k+1}(\mathfrak{m}_\lambda)] = 
\lim_{n\to\infty} n^{-\frac{k+1}{2}} \mathbb{E}_{M^{(n)}_{tn}}[ \Sigma_k ] = 
e^{-kt/m} R_{k+1}(\mathfrak{m}_{\omega_0}) 
\end{equation}
for $k\geqq 2$, and hence \eqref{eq:1-18}. 
This completes the proof of Theorem~\ref{th:1-1}.

\subsection{Proof of Theorem~\ref{th:1-2}}

The verification in the cases of (i) and (ii) goes on similarly to the preceding subsection, proof of Theorem~\ref{th:1-1}, 
by using \eqref{eq:2-8} and \eqref{eq:2-9} in Proposition~\ref{prop:2-2} instead of Proposition~\ref{prop:2-1}. 

Let us consider the case of (iii). 
Similarly to \eqref{eq:2-1-2} and \eqref{eq:2-1-3} in the proof of Theorem~\ref{th:1-1}, 
\eqref{eq:1-21} is derived from \eqref{eq:2-10} in Proposition~\ref{prop:2-2}. 
For simplicity set 
\[ 
g_\alpha (u) = \frac{\sin (\pi\alpha)}{\pi\alpha} \int_0^\infty 
e^{-u (\xi\cos(\pi\alpha/2))^{1/\alpha}} \frac{1}{\xi^2+2\xi\cos (\pi\alpha) +1}\, d\xi, 
\qquad u>0. 
\] 
For a sequence $\{(\mathbb{Y}_n, M^{(n)}_{t\theta_n})\}_{n\in\mathbb{N}}$, 
the same argument to \eqref{eq:2-1-1} yields for $\rho, \sigma\in\mathbb{Y}^\times$ 
\begin{align*}
&\mathbb{E}_{M^{(n)}_{t\theta_n}} [ \Sigma_{\rho\sqcup\sigma}] - 
\mathbb{E}_{M^{(n)}_{t\theta_n}} [ \Sigma_\rho] \, \mathbb{E}_{M^{(n)}_{t\theta_n}} [ \Sigma_\sigma] \\ 
&= f(|\rho|+|\sigma|, n, t\theta_n) \mathbb{E}_{M^{(n)}_0} [ \Sigma_{\rho\sqcup\sigma}] - 
f(|\rho|, n, t\theta_n) f(|\sigma|, n, t\theta_n) \, \mathbb{E}_{M^{(n)}_0} [ \Sigma_\rho] \, 
\mathbb{E}_{M^{(n)}_0} [ \Sigma_\sigma] \\ 
&= \bigl( f(|\rho|+|\sigma|, n, t\theta_n) - g_\alpha (t(|\rho|+|\sigma|)^{1/\alpha})\bigr) \, 
\mathbb{E}_{M^{(n)}_0} [ \Sigma_{\rho\sqcup\sigma}] \\ 
&\quad + \bigl( g_\alpha (t|\rho|^{1/\alpha}) -f(|\rho|, n, t\theta_n)\bigr) f(|\sigma|, n, t\theta_n) 
\mathbb{E}_{M^{(n)}_0} [ \Sigma_\rho] \, \mathbb{E}_{M^{(n)}_0} [ \Sigma_\sigma] \\ 
&\quad + g_\alpha (t|\rho|^{1/\alpha}) \bigl( g_\alpha (t|\sigma|^{1/\alpha}) -f(|\sigma|, n, t\theta_n)\bigr) \, 
\mathbb{E}_{M^{(n)}_0} [ \Sigma_\rho] \, \mathbb{E}_{M^{(n)}_0} [ \Sigma_\sigma] \\ 
&\quad + g_\alpha (t(|\rho|+|\sigma|)^{1/\alpha}) \mathbb{E}_{M^{(n)}_0} [ \Sigma_{\rho\sqcup\sigma}]
- g_\alpha (t|\rho|^{1/\alpha}) g_\alpha (t|\sigma|^{1/\alpha}) \, 
\mathbb{E}_{M^{(n)}_0} [ \Sigma_\rho] \, \mathbb{E}_{M^{(n)}_0} [ \Sigma_\sigma] \\ 
&= \bigl\{ g_\alpha (t(|\rho|+|\sigma|)^{1/\alpha}) - g_\alpha (t|\rho|^{1/\alpha}) g_\alpha (t|\sigma|^{1/\alpha}) \bigr\} 
\mathbb{E}_{M^{(n)}_0} [ \Sigma_\rho] \, \mathbb{E}_{M^{(n)}_0} [ \Sigma_\sigma] \\ 
&\quad + o\bigl( n^{\frac{1}{2}(|\rho|+l(\rho)+|\sigma|+l(\sigma))} \bigr) \qquad 
\text{as} \quad n\to\infty
\end{align*}
by \eqref{eq:2-10} with \eqref{eq:2-2} and \eqref{eq:2-31} for $M^{(n)}_0$. 
Moreover, from \eqref{eq:2-2} and \eqref{eq:2-3} for $M^{(n)}_0$, this equals 
\begin{align}
&= \bigl\{ g_\alpha (t(|\rho|+|\sigma|)^{1/\alpha}) - g_\alpha (t|\rho|^{1/\alpha}) g_\alpha (t|\sigma|^{1/\alpha}) \bigr\} 
\Bigl(\prod_{i=1}^{l(\rho)} r_{\rho_i+1}\Bigr) \Bigl(\prod_{i=1}^{l(\sigma)} r_{\sigma_i+1}\Bigr)  
n^{\frac{1}{2}(|\rho|+l(\rho)+|\sigma|+l(\sigma))} \notag \\ 
&\quad + o\bigl( n^{\frac{1}{2}(|\rho|+l(\rho)+|\sigma|+l(\sigma))} \bigr) \qquad 
\text{as} \quad n\to\infty. 
\label{eq:2-2-1}
\end{align}
If $\omega_0 = \varOmega$ is the initial profile, \eqref{eq:2-2-1} contains only the error term since 
the $j$th free cumulant of $\mathfrak{m}_\varOmega$ vanishes for $j\geqq 3$. 
If $\omega\neq\varOmega$, there exists $k\geqq 2$ such that $r_{k+1} \neq 0$. 
In the case of $\rho = \sigma = (k^m)$ ($m\in\mathbb{N}$) in \eqref{eq:2-2-1}, the main term is 
\[ 
\bigl\{ g_\alpha (t(2km)^{1/\alpha}) - g_\alpha (t(km)^{1/\alpha})^2\bigr\}\, r_{k+1}^{2m} \,
n^{\frac{1}{2}(|\rho|+l(\rho)+|\sigma|+l(\sigma))}.
\] 
As verified below, for any $t>0$, appropriately taken $m$ yields $g_\alpha (t(2km)^{1/\alpha})\neq g_\alpha (t(km)^{1/\alpha})^2$. 
Then the main term does not vanish in \eqref{eq:2-2-1}, which implies that 
$\{(\mathbb{Y}_n, M^{(n)}_{t\theta_n})\}_{n\in\mathbb{N}}$ does not satisfy approximate factorization property. 
Since $g_\alpha (u)$ satisfies 
%
%\begin{equation}\label{eq:2-2-2}
\[ 
g_\alpha (u) \ \sim \ \frac{2}{\pi} \varGamma(\alpha) \sin\frac{\pi\alpha}{2} \, \frac{1}{u^\alpha}
\qquad \text{as} \quad u\to\infty, 
\] 
%\end{equation}
%
we have for $p\in\mathbb{N}$ 
\[ 
\frac{g_\alpha (t(2p)^{1/\alpha})}{g_\alpha (tp^{1/\alpha})^2} \ \sim \ 
\frac{\pi t^\alpha p}{4\varGamma(\alpha) \sin(\pi\alpha/2)} 
\qquad \text{as} \quad p\to\infty. 
\] 
In particular, this is larger than $1$ if $p$ is large enough. 
This completes the proof of Theorem~\ref{th:1-2}.

\section{Technical details}

First we show an inversion formula expressing $f(k, n, s)$ in \eqref{eq:2-7} in term of the characteristic 
function of $\psi$.

\begin{lemma}\label{lem:3-1}
Let $n, k\in\mathbb{N}$ such that $n\geqq k\geqq 2$. 
We have 
\begin{align}
&f(k, n, s) + \frac{k}{2n} \sum_{j=0}^\infty \Bigl( 1- \frac{k}{n}\Bigr)^j \psi^{\ast (j+1)}(\{s\}) 
\notag \\ 
&= \lim_{\epsilon\downarrow 0, r\uparrow \infty} \frac{1}{2\pi i} \int_{\{\epsilon<|\xi|<r\}} 
\frac{e^{-i\xi s}}{1- \varphi(\xi) + (k/n)\varphi(\xi)} \,\frac{\varphi(\xi) -1}{\xi} \,d\xi, \qquad 
s\geqq 0 \label{eq:3-0-1}
\end{align}
where 
\[ 
\varphi(\xi) = \int_0^\infty e^{i\xi u}\psi(du). 
\] 
\end{lemma}
\textbf{Proof} \ 
We compute 
\begin{equation}\label{eq:3-0-2}
\frac{1}{2\pi} \int_{\{\epsilon<|\xi|< r\}} e^{-i\xi s} \Bigl( \int_0^a f(k, n, x) e^{i\xi x} dx \Bigr) d\xi, 
\qquad a>0 
\end{equation}
in two ways. 
Set $f(k, n, s) =0$ for $s<0$ for convenience. 

On one hand, we show 
\begin{equation}\label{eq:3-0-3}
\lim_{\epsilon\downarrow 0, r\uparrow \infty} \lim_{a\uparrow\infty} \text{\eqref{eq:3-0-2}} = 
f(k, n, s)+ \frac{k}{2n} \sum_{j=0}^\infty \Bigl( 1- \frac{k}{n}\Bigr)^j \psi^{\ast (j+1)}(\{s\}), 
\qquad s\in\mathbb{R}.
\end{equation}
(Here we do not care about general conditions for $f(k, n,\,\cdot\,)$ yielding the \lq inversion\rq\  \eqref{eq:3-0-3} 
but use the special form \eqref{eq:2-7} of our $f(k, n, s)$.) 
In 
\begin{equation}\label{eq:3-0-4}
\text{\eqref{eq:3-0-2}} = \frac{1}{\pi} \int_0^a f(k, n, x) \Bigl( \frac{\sin r(x-s)}{x-s} - \frac{\sin \epsilon(x-s)}{x-s}
\Bigr) dx, 
\end{equation}
putting \eqref{eq:2-7} then interchanging the integral and infinite sum, we have
\begin{align*}
&\int_0^a f(k, n, x) \frac{\sin r(x-s)}{x-s} dx \\ 
&= \sum_{j=0}^\infty \Bigl(1-\frac{k}{n}\Bigr)^j 
\int_0^a \Bigl( \int_{[0, x]} \frac{\sin r(x-s)}{x-s} \,\psi\bigl( (x-u, \infty)\bigr) \psi^{\ast j}(du) \Bigr) dx \\ 
&= \sum_{j=0}^\infty \Bigl(1-\frac{k}{n}\Bigr)^j \iint_{[0, \infty)^2} 
\Bigl( \int_0^a \frac{\sin r(x-s)}{x-s} 1_{[0, x]}(u) 1_{(x-u, \infty)}(v)dx \Bigr) \psi^{\ast j}\times\psi (dudv).
\end{align*}
Let $\bigstar$ be the above three-fold integral. 
Then, since 
\begin{align*} 
&\sup_{\alpha < \beta} \Bigl| \int_\alpha^\beta \frac{\sin r(x-s)}{x-s} dx \Bigr| < \infty \qquad 
(\text{also independent of } r, s) \quad \text{and} \\ 
&\int_\alpha^\beta \frac{\sin rx}{x} dx \xrightarrow[\,r\uparrow\infty\,]{} 
\begin{cases}0, & \alpha<\beta<0, \ 0<\alpha<\beta, \\ \pi/2, & \alpha<\beta =0, \ 0=\alpha<\beta, \\ 
\pi, & \alpha<0<\beta \end{cases} 
\end{align*}
hold, the convergence theorem for integral gives 
\begin{align*}
\bigstar &= \iint_{[0, \infty)^2} 
\Bigl( \int_{u\wedge a}^{(u+v)\wedge a} \frac{\sin r(x-s)}{x-s} dx \Bigr) \psi^{\ast j}\times\psi (dudv) \\ 
&\xrightarrow[\,a\uparrow\infty\,]{} \iint_{[0, \infty)^2} \Bigl( \int_u^{u+v} \frac{\sin r(x-s)}{x-s} dx \Bigr) 
\psi^{\ast j}\times\psi (dudv)  \\ 
&= \iint_{[0, \infty)^2} \Bigl( \int_{u-s}^{u+v-s} \frac{\sin ry}{y} dy \Bigr) 
\psi^{\ast j}\times\psi (dudv)  \\ 
&\xrightarrow[\,r\uparrow\infty\,]{} \frac{\pi}{2} \psi^{\ast j}(\{s\}) \psi\bigl( (0,\infty)\bigr) + 
\frac{\pi}{2} \psi^{\ast j}\times\psi(\{u+v =s, v>0\})  \\ 
&\qquad\qquad + \pi \psi^{\ast j}\times\psi (\{ u<s<u+v, v>0\}). 
\end{align*}
The third term is rewritten as 
\[ 
\pi \int_{[0, s]} \psi\bigl( (s-u, \infty)\bigr) \psi^{\ast j}(du) - \pi \psi\bigl( (0,\infty)\bigr) \psi^{\ast j}(\{s\}). 
\] 
After replacing $r$ by $\epsilon$, the part of $a\uparrow\infty$ is the same. 
Then, by $\lim_{\epsilon\downarrow 0}\int_\alpha^\beta \sin\epsilon x /x dx =0$, 
we see from the convergence theorem 
\begin{align*}
&\lim_{\epsilon\downarrow 0, r\uparrow \infty} \lim_{a\uparrow\infty} \text{\eqref{eq:3-0-4}} \\
&= 
\sum_{j=0}^\infty \bigl( 1-\frac{k}{n}\bigr)^j \Bigl\{ \int_{[0, s]} \psi\bigl( (s-u, \infty)\bigr) 
\psi^{\ast j}(du) +\frac{1}{2}\psi^{\ast j}\times\psi(\{u+v=s, v>0\}) \\ 
&\qquad\qquad\qquad\qquad\qquad -\frac{1}{2} \psi^{\ast j}(\{s\}) \psi\bigl((0, \infty)\bigr)\Bigr\} \\ 
&= f(k, n, s) + \frac{1}{2}\sum_{j=0}^\infty \bigl( 1-\frac{k}{n}\bigr)^j 
\Bigl( \int_{(0, s]} \psi^{\ast j}(\{s\}-v) \psi (dv) - \psi^{\ast j}(\{s\}) \psi\bigl((0, \infty)\bigr) \Bigr) \\ 
&= f(k,n,s) + \frac{1}{2}\sum_{j=0}^\infty \bigl( 1-\frac{k}{n}\bigr)^j 
\bigl( \psi^{\ast (j+1)}(\{s\}) - \psi^{\ast j}(\{s\}) \bigr),  
%\\ &= f(k,n,s) + \frac{k}{2n} \sum_{j=0}^\infty \bigl( 1-\frac{k}{n}\bigr)^j \psi^{\ast (j+1)}(\{s\}), 
\end{align*}
which agrees with \eqref{eq:3-0-3}.

On the other hand, the convergence theorem yields 
\begin{equation}\label{eq:3-0-5}
\int_0^a f(k, n, x) e^{i\xi x} dx = \sum_{j=0}^\infty \bigl( 1-\frac{k}{n}\bigr)^j \int_0^a e^{i\xi x} 
\Bigl( \int_{[0, x]} \psi\bigl((x-u, \infty)\bigr) \psi^{\ast j}(du) \Bigr) dx. 
\end{equation}
Let $\bigstar^\prime$ be the two-fold integral in \eqref{eq:3-0-5}. 
We have for $\xi\neq 0$ 
\begin{align}
\bigstar^\prime &= \iint_{[0, \infty)^2} \Bigl( \int_0^a e^{i\xi x} 1_{[0, x]}(u) 1_{(x-u, \infty)}(v) dx \Bigr) 
\psi^{\ast j}\times\psi (dudv) \notag\\ 
&= \iint_{[0, \infty)^2} \Bigl( \int_u^{(u+v)\wedge a} e^{i\xi x} dx \Bigr) 1_{[0, a]}(u) 
\psi^{\ast j}\times\psi (dudv) \notag\\ 
&= \iint_{[0, \infty)^2} \frac{1}{i\xi} \bigl(e^{i\xi ((u+v)\wedge a)} - e^{i\xi u}\bigr) 1_{[0, a]}(u) 
\psi^{\ast j}\times\psi (dudv) \label{eq:3-0-6}\\ 
&\xrightarrow[\,a\uparrow\infty\,]{} 
\iint_{[0, \infty)^2} \frac{1}{i\xi} \bigl(e^{i\xi (u+v)} - e^{i\xi u}\bigr) 
\psi^{\ast j}\times\psi (dudv) 
= \frac{1}{i\xi} \varphi(\xi)^j \bigl( \varphi(\xi) -1\bigr). \notag
%\int_0^a \Bigl( \int_u^a \Bigl( \int_{x-u}^\infty e^{i\xi x} \psi^{\ast j}(u) \psi(v) dv\Bigr)dx\Bigr)du 
%\notag \\ 
%&= \int_0^a \Bigl\{ \int_0^a \Bigl( \int_u^{u+v} e^{i\xi x} \psi^{\ast j}(u) \psi(v) dx \Bigr) dv + 
%\int_a^\infty \Bigl( \int_u^a e^{i\xi x} \psi^{\ast j}(u) \psi(v) dx\Bigr) dv \Bigr\} du 
%\notag \\ 
%&= \frac{1}{i\xi} \Bigl( \int_0^a \psi^{\ast j}(u) e^{i\xi u} du \Bigr) 
%\Bigl( \int_0^a \psi(v) (e^{i\xi v} -1) dv\Bigr) 
%\notag \\ 
%&\quad + \frac{1}{i\xi} \Bigl( \int_a^\infty \psi(v)dv\Bigr) \Bigl\{ e^{ia\xi} 
%\Bigl( \int_0^a \psi^{\ast j}(u)du\Bigr) - \Bigl( \int_0^a \psi^{\ast j}(u)e^{i\xi u} du\Bigr) \Bigr\} 
%\label{eq:3-0-6} \\ 
%&\xrightarrow[\,a\uparrow\infty\,]{} \frac{1}{i\xi} \varphi(\xi)^j \bigl( \varphi(\xi) -1\bigr). 
%\notag 
\end{align}
Hence, from the convergence theorem, 
\[ 
\text{\eqref{eq:3-0-5}} \xrightarrow[\,a\uparrow\infty\,]{} \sum_{j=0}^\infty\bigl(1-\frac{k}{n}\bigr)^j 
\frac{1}{i\xi}\varphi(\xi)^j \bigl(\varphi(\xi)-1\bigr) = 
\frac{\varphi(\xi) -1}{i\xi} \,\frac{1}{1- (1- (k/n))\varphi(\xi)}.
\] 
Since we see from the expression of \eqref{eq:3-0-6} that the left hand side of \eqref{eq:3-0-5} is bounded 
jointly with respect to $\{\epsilon<|\xi|<r\}$ and $a>0$, 
we have 
\[ 
\text{\eqref{eq:3-0-2}} \ \xrightarrow[\,a\uparrow\infty\,] \ \frac{1}{2\pi} \int_{\{\epsilon<|\xi|<r\}} 
e^{-i\xi s}\, \frac{\varphi(\xi) -1}{i\xi} \,\frac{1}{1- (1- (k/n))\varphi(\xi)} d\xi
\] 
by the convergence theorem for integral. 
Combined with the former half, this completes the proof of \eqref{eq:3-0-1}. 
\hfill $\blacksquare$

\medskip

Before entering into the proof of Proposition~\ref{prop:2-1}, let us note the case 
where the pausing time obeys an exponential distribution: 
\[ 
\psi(dx) = \frac{1}{m} e^{-x/m}\, 1_{[0, \infty)}(x)dx, \qquad 
\varphi(\xi) = \frac{1}{1-im\xi}.
\] 
Then Lemma~\ref{lem:3-1} gives (with residue calculus) 
\begin{equation}\label{eq:3-0-7}
f(k, n, tn) = 
\lim_{r\uparrow\infty} \frac{m}{2\pi} \int_{\{|\xi|<r\}} \frac{e^{-itn\xi}}{\frac{k}{n} -im\xi} d\xi = 
\lim_{r\uparrow\infty}\frac{1}{2\pi i} \int_{\{|\xi|<nr\}} \frac{e^{-it\xi}}{-i\frac{k}{m}-\xi} d\xi = 
e^{-\frac{kt}{m}}.
\end{equation}

\subsection{Proof of Proposition~\ref{prop:2-1}}

Under the assumptions of Proposition~\ref{prop:2-1}, $\psi$ is a continuous distribution. 
In fact, the integrability of \eqref{eq:1-15} and uniform continuity of $\varphi$ yield 
$\lim_{\xi\to\pm\infty} |\varphi(\xi)| =0$, from which continuity of $\psi$ follows 
(see \cite[\S2.2]{Lu60}). 

Now the atomic parts do not appear in \eqref{eq:3-0-1} of Lemma~\ref{lem:3-1}. 
Since the integrand of \eqref{eq:3-0-1} does not have singularity as $\epsilon\downarrow 0$ by the 
differentiability of $\varphi$ at $0$, we have 
\begin{equation}\label{eq:3-1-1}
f(k, n, tn) = \lim_{r\uparrow\infty} \frac{-1}{2\pi i} \int_{\{|\xi|<r\}} \frac{e^{-itn\xi}}{\xi}\, 
\frac{1-\varphi(\xi)}{1-\varphi(\xi)+ (k/n) \varphi(\xi)} d\xi.
\end{equation}
We divide the integral of \eqref{eq:3-1-1} into the following four pieces where $\delta >0$ is specified a bit later: 
\begin{align}
&\frac{-1}{2\pi i} \int_{\{|\xi|<r\}} \frac{e^{-itn\xi}}{\xi} 
\frac{1-\varphi(\xi)}{1-\varphi(\xi)+ (k/n) \varphi(\xi)} d\xi \notag\\ 
&= \frac{-1}{2\pi i} \Bigl\{\int_{\{|\xi|\leqq \frac{\delta}{n}\}} \frac{e^{-itn\xi}}{\xi} 
\frac{1-\varphi(\xi)}{1-\varphi(\xi)+ \frac{k}{n} \varphi(\xi)} d\xi + 
\int_{\{\frac{\delta}{n}<|\xi|<r\}} \frac{e^{-itn\xi}}{\xi} 
\Bigl( 1- \frac{\frac{k}{n}\varphi(\xi)}{1-\varphi(\xi)+\frac{k}{n}\varphi(\xi)} \Bigr) d\xi \Bigr\} \notag\\ 
&= \frac{-1}{2\pi i} \int_{\{|\xi|\leqq \frac{\delta}{n}\}} \frac{e^{-itn\xi}}{\xi} 
\frac{1-\varphi(\xi)}{1-\varphi(\xi)+ \frac{k}{n} \varphi(\xi)} d\xi \ 
- \frac{1}{2\pi i} \int_{\{\frac{\delta}{n}<|\xi|<r\}} \frac{e^{-itn\xi}}{\xi} d\xi \notag\\ 
&\quad + \frac{k}{2\pi i} \int_{\{\frac{\delta}{n}<|\xi|<\delta\}} 
\frac{e^{-itn\xi}}{\xi} \frac{\varphi(\xi)}{n(1-\varphi(\xi)) +k\varphi(\xi)} d\xi \notag\\ 
&\quad + \frac{k}{2\pi i} \int_{\{\delta\leqq |\xi|<r\}} 
\frac{e^{-itn\xi}}{\xi} \frac{\varphi(\xi)}{n(1-\varphi(\xi)) +k\varphi(\xi)} d\xi \notag\\ 
&= \text{(I)} + \text{(II)} + \text{(III)} + \text{(IV)}. 
\label{eq:3-1-2}
\end{align}

First we look at (IV) in \eqref{eq:3-1-2}. 
We have 
\[ 
|n(1-\varphi(\xi))+k\varphi(\xi)| \geqq n|1-\varphi(\xi)|-k|\varphi(\xi)| \geqq 
n \bigl( 1- \sup_{\delta\leqq |\xi|} |\varphi(\xi)| \bigr) -k, 
\] 
in which $1- \sup_{\delta\leqq |\xi|} |\varphi(\xi)| >0$ holds for any $\delta >0$. 
In fact, since we saw $\lim_{\xi\to\pm\infty}\varphi(\xi) =0$, take $\delta^\prime >0$ such that 
$\sup_{|\xi|>\delta^\prime} |\varphi(\xi)| \leqq 1/2$. 
If $\delta<\delta^\prime$ and $\sup_{\delta\leqq |\xi|\leqq\delta^\prime} |\varphi(\xi)| =1$, 
we choose a sequence $\{\xi_n\} \subset \{\delta\leqq |\xi|\leqq\delta^\prime\}$ such that 
$\lim_{n\to\infty} |\varphi(\xi_n)|=1$. 
There exists $\xi_0 \in\{\delta\leqq |\xi|\leqq\delta^\prime\}$ such that $|\varphi(\xi_0)| =1$ by the compactness, 
namely $\int_{\mathbb{R}} e^{i\xi_0 x} \psi(dx) = e^{ia}$ for some $a\in\mathbb{R}$. 
We have, however, 
\[ 
\int_{\mathbb{R}} |e^{i\xi_0 x} -e^{ia}|^2 \psi(dx) = 1-1-1+1 =0, 
\] 
contradicting continuity of $\psi$. 
We thus obtain 
\begin{equation}\label{eq:3-1-3}
\bigl| \text{(IV) in \eqref{eq:3-1-2}}\bigr| \leqq \frac{k}{2\pi \{n(1- \sup_{\delta\leqq |\xi|} |\varphi(\xi)|)-k\}} \,
\int_{\{\delta\leqq |\xi|\}} \Bigl| \frac{\varphi(\xi)}{\xi}\Bigr| d\xi 
\end{equation}
for any $\delta >0$ (with the upper bound independent of $r$). 
We note 
\[
\text{(II) in \eqref{eq:3-1-2}} = - \frac{1}{2\pi i} \int_{\{\delta<|\xi|<nr\}} \frac{e^{-it\xi}}{\xi} d\xi 
\] 
converges as $r\uparrow \infty$. 

Let us compute 
\begin{equation}\label{eq:3-1-4}
\text{(III) in \eqref{eq:3-1-2}} = \frac{k}{2\pi i} \int_{\mathbb{R}} 1_{\{\delta<|\xi|<n\delta\}}(\xi) 
\frac{e^{-it\xi}}{\xi} \frac{\varphi(\xi/n)}{n(1-\varphi(\xi/n))+k\varphi(\xi/n)} d\xi.
\end{equation}
Noting $\varphi^\prime (0) = im$, we have for any $\delta >0$ and any $\xi\in\mathbb{R}$ such that $\delta<|\xi|$ 
\[ 
\text{integrand of \eqref{eq:3-1-4}} \, \xrightarrow[n\to\infty]{}\, \frac{e^{-it\xi}}{\xi} \frac{1}{-im\xi +k}.
\] 
To verify uniform integrability of the integrand, $|\xi|<n\delta$ yields 
\begin{align*}
&\bigl| n\bigl( 1-\varphi(\frac{\xi}{n}) \bigr) + k\varphi(\frac{\xi}{n}) \bigr| = 
\bigl| n\bigl( 1-\varphi(\frac{\xi}{n}) \bigr) + \varphi^\prime (0)\xi -im\xi + k\varphi(\frac{\xi}{n}) -k+k\bigr| \\ 
&\geqq |-im\xi+k| - \bigl| n\bigl( 1-\varphi(\frac{\xi}{n}) \bigr) + \varphi^\prime (0)\xi 
\bigr| - k\bigl|\varphi(\frac{\xi}{n})-1\bigr| \\ 
&= \sqrt{m^2 \xi^2 +k^2} - \Bigl| \frac{\varphi(\xi/n)-1}{\xi/n} - \varphi^\prime (0)\Bigr| |\xi| 
- k\bigl|\varphi(\frac{\xi}{n})-1\bigr| \\ 
&\geqq \frac{1}{\sqrt{2}}(m |\xi| +k) - A_\delta |\xi| - B_\delta k 
\end{align*}
by setting 
\[ 
A_\delta = \sup_{|\eta|<\delta} \Bigl| \frac{\varphi(\eta)-1}{\eta} -\varphi^\prime (0)\Bigr|, \qquad 
B_\delta = \sup_{|\eta|<\delta} |\varphi(\eta)-1|.
\] 
For any $\delta >0$ such that 
\begin{equation}\label{eq:3-1-5}
A_\delta < m/\sqrt{2}, \qquad B_\delta < 1/\sqrt{2},  
\end{equation}
we have 
\[ 
\bigl|\text{integrand of \eqref{eq:3-1-4}}\bigr| \leqq 1_{\{\delta<|\xi|\}}(\xi) \frac{1}{|\xi|} \, 
\frac{1}{((m/\sqrt{2})- A_\delta) |\xi| + ((1/\sqrt{2})-B_\delta) k}, 
\] 
the right hand side being an integrable function in $\xi$. 
If $\delta >0$ satisfies \eqref{eq:3-1-5}, then 
\begin{align}
\text{\eqref{eq:3-1-4}} &\xrightarrow[n\to\infty]{} 
\frac{k}{2\pi i} \int_{\mathbb{R}} 1_{\{\delta<|\xi|\}}(\xi) \frac{e^{-it\xi}}{\xi (-im\xi +k)} d\xi 
= \frac{1}{2\pi i} \int_{\{\delta<|\xi|\}} e^{-it\xi} \Bigl( \frac{1}{\xi} + \frac{im}{k-im\xi}\Bigr) d\xi \notag\\ 
&= \frac{1}{2\pi i} \lim_{R\to\infty} \int_{\{\delta<|\xi|<R\}} \frac{e^{-it\xi}}{\xi} d\xi + 
\frac{1}{2\pi} \lim_{R\to\infty} \int_{\{\delta<|\xi|<R\}} \frac{e^{-it\xi}}{(k/m)-i\xi} d\xi \notag\\ 
&= \frac{1}{2\pi i} \lim_{R\to\infty} \int_{\{\delta<|\xi|<R\}} \frac{e^{-it\xi}}{\xi} d\xi + 
e^{-kt/m} + O(\delta) \label{eq:3-1-6}
\end{align}
(recall the computation of \eqref{eq:3-0-7}). 

Let us apply these estimates to \eqref{eq:3-1-1} and \eqref{eq:3-1-2}: 
\[ 
f(k,n, tn) = \text{(I)} + \text{(III)} + \lim_{r\uparrow \infty} \text{(II)} + \lim_{r\uparrow\infty} \text{(IV)}.
\] 
For arbitrarily given $\epsilon >0$, take $\delta >0$ satisfying $\text{(I)}<\epsilon$, 
$O(\delta)$ in \eqref{eq:3-1-6} $<\epsilon$ and \eqref{eq:3-1-5}. 
We then take sufficiently large $N\in\mathbb{N}$ such that, if $n\geqq N$, 
\begin{align*}
&\bigl| \lim_{r\to\infty} \text{(IV)} \bigr| < \epsilon \qquad (\text{by \eqref{eq:3-1-3}}) \\ 
&\text{(III)} = \frac{1}{2\pi i} \lim_{R\to\infty} \int_{\{\delta<|\xi|<R\}} \frac{e^{-it\xi}}{\xi} d\xi + 
e^{-kt/m} + \; (\text{term of } |~\cdot~|\leqq 2\epsilon).
\end{align*}
The first term agrees with $-\lim_{r\uparrow\infty}\text{(II)}$. 
Consequently, we obtain 
\[ 
\limsup_{n\to\infty} \bigl| f(k, n, tn) - e^{-kt/m} \bigr| \leqq 4 \epsilon
\] 
for any $\epsilon >0$. 
This completes the proof of Proposition~\ref{prop:2-1}.

\subsection{Proof of Proposition~\ref{prop:2-2}}

\begin{figure}[hbt]
\begin{center}
{\unitlength 0.1in%
\begin{picture}(20.3000,20.2700)(-0.3000,-24.0000)%
% LINE 2 0 3 0 Black White  
% 4 600 2400 600 1400 600 1400 400 1400
% 
\special{pn 8}%
\special{pa 600 2400}%
\special{pa 600 1400}%
\special{fp}%
\special{pa 600 1400}%
\special{pa 400 1400}%
\special{fp}%
% VECTOR 2 0 3 0 Black White  
% 4 600 1400 600 400 600 1400 2000 1400
% 
\special{pn 8}%
\special{pa 600 1400}%
\special{pa 600 400}%
\special{fp}%
\special{sh 1}%
\special{pa 600 400}%
\special{pa 580 467}%
\special{pa 600 453}%
\special{pa 620 467}%
\special{pa 600 400}%
\special{fp}%
\special{pa 600 1400}%
\special{pa 2000 1400}%
\special{fp}%
\special{sh 1}%
\special{pa 2000 1400}%
\special{pa 1933 1380}%
\special{pa 1947 1400}%
\special{pa 1933 1420}%
\special{pa 2000 1400}%
\special{fp}%
% CIRCLE 0 0 3 0 Black White  
% 4 600 1400 1800 1400 1670 900 1330 530
% 
\special{pn 20}%
\special{ar 600 1400 1200 1200 5.4105094 5.8460466}%
% CIRCLE 0 0 3 0 Black White  
% 4 600 1400 1800 1400 1380 2330 1690 1910
% 
\special{pn 20}%
\special{ar 600 1400 1200 1200 0.4376312 0.8728935}%
% CIRCLE 0 0 3 0 Black White  
% 4 600 1400 600 1400 600 1400 600 1400
% 
\special{pn 20}%
\special{ar 600 1400 0 0 0.0000000 6.2831853}%
% CIRCLE 0 0 3 0 Black White  
% 4 600 1400 800 1400 760 1340 730 1250
% 
\special{pn 20}%
\special{ar 600 1400 200 200 5.4264797 5.9244146}%
% CIRCLE 0 0 3 0 Black White  
% 4 600 1400 800 1400 710 1530 770 1480
% 
\special{pn 20}%
\special{ar 600 1400 200 200 0.4398426 0.8685394}%
% LINE 0 0 3 0 Black White  
% 2 730 1230 1360 490
% 
\special{pn 20}%
\special{pa 730 1230}%
\special{pa 1360 490}%
\special{fp}%
% LINE 0 0 3 0 Black White  
% 2 800 1340 1670 890
% 
\special{pn 20}%
\special{pa 800 1340}%
\special{pa 1670 890}%
\special{fp}%
% LINE 0 0 3 0 Black White  
% 2 730 1570 1350 2330
% 
\special{pn 20}%
\special{pa 730 1570}%
\special{pa 1350 2330}%
\special{fp}%
% LINE 0 0 3 0 Black White  
% 2 800 1470 1690 1900
% 
\special{pn 20}%
\special{pa 800 1470}%
\special{pa 1690 1900}%
\special{fp}%
% STR 2 0 3 0 Black White  
% 4 1370 1600 1370 1700 2 0 0 0
% $L_-$
\put(13.7000,-17.0000){\makebox(0,0)[lb]{$L_-$}}%
% STR 2 0 3 0 Black White  
% 4 1130 1960 1130 2060 4 0 0 0
% $L_-^\prime$
\put(11.3000,-20.6000){\makebox(0,0)[rt]{$L_-^\prime$}}%
% STR 2 0 3 0 Black White  
% 4 1600 1990 1600 2090 1 0 0 0
% $C_-$
\put(16.0000,-20.9000){\makebox(0,0)[lt]{$C_-$}}%
% STR 2 0 3 0 Black White  
% 4 740 1510 740 1610 3 0 0 0
% $C_-^\prime$
\put(7.4000,-16.1000){\makebox(0,0)[rb]{$C_-^\prime$}}%
% STR 2 0 3 0 Black White  
% 4 1370 980 1370 1080 1 0 0 0
% $L_+$
\put(13.7000,-10.8000){\makebox(0,0)[lt]{$L_+$}}%
% STR 2 0 3 0 Black White  
% 4 1050 720 1050 820 3 0 0 0
% $L_+^\prime$
\put(10.5000,-8.2000){\makebox(0,0)[rb]{$L_+^\prime$}}%
% STR 2 0 3 0 Black White  
% 4 1550 570 1550 670 2 0 0 0
% $C_+$
\put(15.5000,-6.7000){\makebox(0,0)[lb]{$C_+$}}%
% STR 2 0 3 0 Black White  
% 4 700 1120 700 1220 4 0 0 0
% $C_+^\prime$
\put(7.0000,-12.2000){\makebox(0,0)[rt]{$C_+^\prime$}}%
\end{picture}}%
\caption{Contours}
\label{fig:1}
\end{center}
\end{figure}
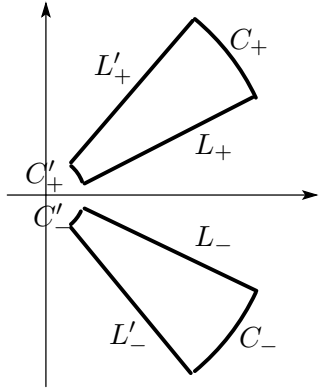

Putting \eqref{eq:1-20}, the expression of the characteristic function, into \eqref{eq:3-0-1} of Lemma~\ref{lem:3-1} 
and noting (absolute) continuity of $\psi$, we have 
\begin{align}
f(k, n, s) &= \lim_{\epsilon\downarrow 0, \, r\uparrow \infty} \frac{1}{2\pi i} \int_{\{\epsilon<|\xi|<r\}} 
\frac{e^{-is\xi}}{\xi} \, 
\frac{e^{-|\xi|^\alpha (1-i\tan(\pi\alpha/2) \mathrm{sgn}\xi)} -1}
{1- (1-\frac{k}{n}) e^{-|\xi|^\alpha (1-i\tan(\pi\alpha/2) \mathrm{sgn}\xi)}} d\xi \notag\\ 
&= \lim_{\epsilon\downarrow 0, \, r\uparrow \infty} \frac{1}{2\pi\alpha i} \Bigl\{ \int_{\epsilon^\alpha}^{r^\alpha} 
\frac{e^{-is\eta^{1/\alpha}}}{\eta} \, 
\frac{e^{-\eta (1-i\tan(\pi\alpha/2))}-1}{1-(1-\frac{k}{n}) e^{-\eta (1-i\tan(\pi\alpha/2))}} d\eta \notag\\ 
&\qquad\qquad - \int_{\epsilon^\alpha}^{r^\alpha} 
\frac{e^{is\eta^{1/\alpha}}}{\eta} \, 
\frac{e^{-\eta (1+i\tan(\pi\alpha/2))}-1}{1-(1-\frac{k}{n}) e^{-\eta (1+i\tan(\pi\alpha/2))}} d\eta \Bigr\}. 
\label{eq:3-2-1}
\end{align}
Let us refer to the two integrals in \eqref{eq:3-2-1} as $\bigstar_-$ and $\bigstar_+$ respectively. 
Take line segments $L_\mp, L_\mp^\prime$ and arcs $C_\mp, C_\mp^\prime$ as in Figure~\ref{fig:1}: 
\begin{align*}
&L_\mp: \epsilon^\alpha \bigl(1\mp i\tan\frac{\pi\alpha}{2}\bigr) \,\longrightarrow\, 
r^\alpha \bigl(1 \mp i\tan\frac{\pi\alpha}{2}\bigr), \\  
&L_\mp^\prime: \epsilon^\alpha \bigl(1\mp i\tan\frac{\pi\alpha}{2}\bigr) e^{\mp i\pi\alpha/2} \,\longrightarrow\, 
r^\alpha \bigl(1\mp i\tan\frac{\pi\alpha}{2}\bigr) e^{\mp i\pi\alpha/2}, \\ 
&C_\mp: r^\alpha \bigl(1\mp i\tan\frac{\pi\alpha}{2}\bigr) \,\longrightarrow\, 
r^\alpha \bigl(1\mp i\tan\frac{\pi\alpha}{2}\bigr) e^{\mp i\pi\alpha/2}, \\ 
&C_\mp^\prime: \epsilon^\alpha \bigl(1\mp i\tan\frac{\pi\alpha}{2}\bigr) \,\longrightarrow\, 
\epsilon^\alpha \bigl(1\mp i\tan\frac{\pi\alpha}{2}\bigr) e^{\mp i\pi\alpha/2}.
\end{align*}
%
%\begin{figure}[hbt]
%\begin{center}
%\input{contour-fig.tex}
%%\vspace{-5mm}
%\caption{Contours}
%\label{fig:1}
%\end{center}
%\end{figure}
%

\noindent 
We have 
\[ 
\bigstar_\mp = \int_{L_\mp} F_\mp (z) dz \qquad \text{where} \qquad 
F_\mp (z) = \frac{e^{\mp is(z/(1\mp i\tan(\pi\alpha/2)))^{1/\alpha}}}{1-(1-\frac{k}{n})e^{-z}}
\frac{e^{-z}-1}{z}.
\] 
Noting 
\[ 
\lim_{r\uparrow\infty} \int_{C_\mp} F_\mp (z) dz =0, \qquad 
\lim_{\epsilon\downarrow 0} \int_{C_\mp^\prime} F_\mp (z) dz =0,
\] 
we have 
\begin{align*}
\lim_{\epsilon\downarrow 0, \, r\uparrow\infty} \bigstar_\mp &= 
\lim_{\epsilon\downarrow 0, \, r\uparrow\infty} \int_{L_\mp^\prime} F_\mp (z) dz \\ 
&= \int_0^\infty \frac{e^{-s x^{1/\alpha}}}{1-(1-\frac{k}{n}) e^{-x (1\mp i\tan(\pi\alpha/2))e^{\mp i\pi\alpha/2}}} \, 
\frac{e^{-x (1\mp i\tan(\pi\alpha/2))e^{\mp i\pi\alpha/2}} -1}{x} \, dx, 
\end{align*}
the last integrals converging absolutely near $0$ and $\infty$. 
Putting this into \eqref{eq:3-2-1}, we get an easier integral expression 
\begin{align*}
f(k, n, s) = \frac{1}{2\pi\alpha i} \int_0^\infty 
\frac{e^{-s x^{1/\alpha}}}{x} &\Bigl\{ 
\frac{e^{-x (1-i\tan(\pi\alpha/2)) e^{-i\pi\alpha/2}}-1}{1-(1-\frac{k}{n}) e^{-x (1-i\tan(\pi\alpha/2)) e^{-i\pi\alpha/2}}} 
\notag\\ 
&- \frac{e^{-x (1+i\tan(\pi\alpha/2)) e^{i\pi\alpha/2}}-1}{1-(1-\frac{k}{n}) e^{-x (1+i\tan(\pi\alpha/2)) e^{i\pi\alpha/2}}} 
\Bigr\} dx,   
%\label{eq:3-2-2}
\end{align*}
furthermore, after a bit of computation 
\begin{equation}\label{eq:3-2-3}
= \frac{k}{\pi\alpha n} \int_0^\infty 
\frac{e^{-sx^{\frac{1}{\alpha}}}\, e^{-x\frac{\cos(\pi\alpha)}{\cos(\pi\alpha/2)}}\sin (2x\sin\frac{\pi\alpha}{2})}
{1- 2(1-\frac{k}{n}) e^{-x\frac{\cos(\pi\alpha)}{\cos(\pi\alpha/2)}}\cos (2x\sin\frac{\pi\alpha}{2}) + 
(1-\frac{k}{n})^2 e^{-2x\frac{\cos(\pi\alpha)}{\cos(\pi\alpha/2)}}} \,\frac{dx}{x}.
%\frac{\sin (2\frac{y}{n}\sin \frac{\pi\alpha}{2})}{\frac{y}{n}} \notag\\ 
%&\qquad\quad \times \Bigl\{ n^2\{e^{2\frac{y}{n}\cos(\pi\alpha)/ \cos\frac{\pi\alpha}{2}} +1 
%-2 e^{\frac{y}{n}\cos(\pi\alpha)/ \cos\frac{\pi\alpha}{2}} \cos (2\frac{y}{n}\sin\frac{\pi\alpha}{2})\} 
%\notag\\ 
%&\qquad\qquad +
%2nk\{ e^{\frac{y}{n}\cos(\pi\alpha)/ \cos\frac{\pi\alpha}{2}} \cos (2\frac{y}{n}\sin\frac{\pi\alpha}{2}) -1\} 
%+k^2 \Bigr\}^{-1} dy.
\end{equation}

\begin{lemma}\label{lem:3-2}
If $\alpha \neq 1/2$, then 
\begin{equation}\label{eq:3-2-4}
\lim_{n\to\infty} \frac{1}{n} \int_b^\infty 
\frac{e^{-x\frac{\cos(\pi\alpha)}{\cos(\pi\alpha/2)}}\, dx}
{\bigl| 1- 2(1-\frac{k}{n}) e^{-x\frac{\cos(\pi\alpha)}{\cos(\pi\alpha/2)}}\cos (2x\sin\frac{\pi\alpha}{2}) + 
(1-\frac{k}{n})^2 e^{-2x\frac{\cos(\pi\alpha)}{\cos(\pi\alpha/2)}}\bigr|} =0
\end{equation}
%
%\begin{multline}\label{eq:3-2-4}
%\lim_{n\to\infty} \int_{bn}^\infty e^{\frac{y}{n}\cos(\pi\alpha)/ \cos\frac{\pi\alpha}{2}} 
%\Bigl| n^2\{e^{2\frac{y}{n}\cos(\pi\alpha)/ \cos\frac{\pi\alpha}{2}} +1 
%-2 e^{\frac{y}{n}\cos(\pi\alpha)/ \cos\frac{\pi\alpha}{2}} \cos (2\frac{y}{n}\sin\frac{\pi\alpha}{2})\} \\ 
%+2nk\{ e^{\frac{y}{n}\cos(\pi\alpha)/ \cos\frac{\pi\alpha}{2}} \cos (2\frac{y}{n}\sin\frac{\pi\alpha}{2}) -1\} 
%+k^2 \Bigr|^{-1}
%=0
%\end{multline}
%% 
holds for any $b>0$.
\end{lemma}
\textbf{Proof} \ 
First let $0<\alpha <1/2$, hence $\cos(\pi\alpha) /\cos(\pi\alpha/2) >0$. 
Then, 
\begin{align*}
|\text{denominator of the integrand}| &\geqq 1+ (1-\frac{k}{n})^2 e^{-2x\frac{\cos(\pi\alpha)}{\cos(\pi\alpha/2)}} 
- 2(1-\frac{k}{n}) e^{-x\frac{\cos(\pi\alpha)}{\cos(\pi\alpha/2)}} \\ 
&= \Bigl( 1- (1-\frac{k}{n}) e^{-x\frac{\cos(\pi\alpha)}{\cos(\pi\alpha/2)}}\Bigr)^2 \geqq 
\Bigl( 1- e^{-b\frac{\cos(\pi\alpha)}{\cos(\pi\alpha/2)}}\Bigr)^2 >0.
\end{align*}
This yields \eqref{eq:3-2-4}.

Secondly let $1/2<\alpha <1$, hence $\cos(\pi\alpha) /\cos(\pi\alpha/2) <0$. 
Then, the integrand equals 
\[ 
\frac{e^{x\frac{\cos(\pi\alpha)}{\cos(\pi\alpha/2)}}}
{\bigl| e^{2x\frac{\cos(\pi\alpha)}{\cos(\pi\alpha/2)}} - 
2(1-\frac{k}{n}) e^{x\frac{\cos(\pi\alpha)}{\cos(\pi\alpha/2)}}\cos (2x\sin\frac{\pi\alpha}{2}) + 
(1-\frac{k}{n})^2 \bigr|}, 
\] 
whose denominator is bounded below by 
\begin{align*}
&e^{2x\frac{\cos(\pi\alpha)}{\cos(\pi\alpha/2)}} + (1-\frac{k}{n})^2 - 
2(1-\frac{k}{n}) e^{x\frac{\cos(\pi\alpha)}{\cos(\pi\alpha/2)}} \\ 
&= \Bigl(1-\frac{k}{n}- e^{x\frac{\cos(\pi\alpha)}{\cos(\pi\alpha/2)}} \Bigr)^2 \geqq 
\Bigl(1-\frac{k}{n}- e^{b\frac{\cos(\pi\alpha)}{\cos(\pi\alpha/2)}} \Bigr)^2 >0
\end{align*}
if $n$ is sufficiently large. 
This yields \eqref{eq:3-2-4}, too.
\hfill $\blacksquare$

\bigskip

Setting $s = t\theta_n$ for $t>0$ and $\beta_n = \theta_n/n^{1/\alpha}$ in \eqref{eq:3-2-3}, 
we seek the limit of 
\begin{multline}\label{eq:3-2-5}
f(k, n, t\theta_n) = 
\frac{k}{\pi\alpha n} \int_0^\infty \frac{e^{-t\beta_n (nx)^{\frac{1}{\alpha}}}}{x} \\
\frac{e^{-x\frac{\cos(\pi\alpha)}{\cos(\pi\alpha/2)}}\sin (2x\sin\frac{\pi\alpha}{2})}
{1- 2(1-\frac{k}{n}) e^{-x\frac{\cos(\pi\alpha)}{\cos(\pi\alpha/2)}}\cos (2x\sin\frac{\pi\alpha}{2}) + 
(1-\frac{k}{n})^2 e^{-2x\frac{\cos(\pi\alpha)}{\cos(\pi\alpha/2)}}} dx.
\end{multline}

First we consider the case of $\alpha\neq 1/2$. 
Since we know the integration on $(ak, \infty)$ tends to $0$ as $n\to\infty$ by Lemma~\ref{lem:3-2}, 
let us compute 
\begin{equation}\label{eq:3-2-6}
\frac{k}{\pi\alpha n} \int_0^{ak}
\frac{e^{-t\beta_n (nx)^{\frac{1}{\alpha}}}\, e^{-x\frac{\cos(\pi\alpha)}{\cos(\pi\alpha/2)}}\sin (2x\sin\frac{\pi\alpha}{2})}
{1- 2(1-\frac{k}{n}) e^{-x\frac{\cos(\pi\alpha)}{\cos(\pi\alpha/2)}}\cos (2x\sin\frac{\pi\alpha}{2}) + 
(1-\frac{k}{n})^2 e^{-2x\frac{\cos(\pi\alpha)}{\cos(\pi\alpha/2)}}} \frac{dx}{x}
\end{equation}
where $a>0$ is specified later. 
We rewrite \eqref{eq:3-2-6} as 
\begin{align}
\frac{k^2}{\pi\alpha} \int_0^{an} &e^{-t\beta_n (ky)^{\frac{1}{\alpha}}}
e^{\frac{k}{n}y\frac{\cos(\pi\alpha)}{\cos(\pi\alpha/2)}}\frac{\sin (2\frac{k}{n}y\sin\frac{\pi\alpha}{2})}{\frac{k}{n}y} 
\notag\\
&\Bigl\{n^2\bigl\{e^{2\frac{k}{n}y\frac{\cos(\pi\alpha)}{\cos(\pi\alpha/2)}} +1 
- 2e^{\frac{k}{n}y\frac{\cos(\pi\alpha)}{\cos(\pi\alpha/2)}} \cos (2\frac{k}{n}y\sin\frac{\pi\alpha}{2})\bigr\} 
\notag\\ 
&\qquad + 
2nk\bigl\{e^{\frac{k}{n}y\frac{\cos(\pi\alpha)}{\cos(\pi\alpha/2)}} \cos (2\frac{k}{n}y\sin\frac{\pi\alpha}{2}) -1\bigr\} 
+k^2\Bigr\}^{-1} dy
\label{eq:3-2-7}
\end{align}
and note that
\begin{align}
&\lim_{x\to 0} \frac{1}{x} \Bigl\{ e^{kx\frac{\cos(\pi\alpha)}{\cos(\pi\alpha/2)}} \cos (2kx\sin\frac{\pi\alpha}{2}) -1 
\Bigr\} = \frac{k\cos(\pi\alpha)}{\cos(\pi\alpha/2)}, \notag\\ 
&\lim_{x\to 0} \frac{1}{x^2} \Bigl\{ e^{2\frac{k}{n}y\frac{\cos(\pi\alpha)}{\cos(\pi\alpha/2)}} +1 
- 2e^{\frac{k}{n}y\frac{\cos(\pi\alpha)}{\cos(\pi\alpha/2)}} \cos (2\frac{k}{n}y\sin\frac{\pi\alpha}{2})\Bigr\} 
= \frac{k^2}{\cos^2(\pi\alpha/2)}.
\label{eq:3-2-8}
\end{align}
For any $\epsilon >0$ there exists $a>0$ such that $0<x\leqq a$ implies 
\begin{align*}
&\Bigl| \frac{1}{x} \Bigl\{ e^{kx\frac{\cos(\pi\alpha)}{\cos(\pi\alpha/2)}} \cos (2kx\sin\frac{\pi\alpha}{2}) -1 
\Bigr\} - \frac{k\cos(\pi\alpha)}{\cos(\pi\alpha/2)} \Bigr| \leqq \epsilon, \notag\\ 
&\Bigl| \frac{1}{x^2} \Bigl\{ e^{2\frac{k}{n}y\frac{\cos(\pi\alpha)}{\cos(\pi\alpha/2)}} +1 
- 2e^{\frac{k}{n}y\frac{\cos(\pi\alpha)}{\cos(\pi\alpha/2)}} \cos (2\frac{k}{n}y\sin\frac{\pi\alpha}{2})\Bigr\} 
- \frac{k^2}{\cos^2(\pi\alpha/2)} \Bigr| \leqq \epsilon.
\end{align*}
Then, $\epsilon >0$ being taken smaller than $k^2 /\cos^2 (\pi\alpha /2)$, it holds in \eqref{eq:3-2-7} 
\begin{align*}
\Bigl| \Bigl\{ ~\cdot~ \Bigr\} \Bigr| &\geqq \Bigl( \frac{k^2}{\cos^2 (\pi\alpha /2)} - \epsilon \Bigr) y^2 
+ k^2 - 2ky \Bigl( k \bigl| \frac{\cos(\pi\alpha)}{\cos(\pi\alpha/2)} \bigr| +\epsilon\Bigr) \\ 
&= k^2 \Bigl\{ \Bigl( \frac{1}{\cos^2 (\pi\alpha /2)} - \frac{\epsilon}{k^2} \Bigr) y^2 
-2 \Bigl( \bigl| \frac{\cos(\pi\alpha)}{\cos(\pi\alpha/2)} \bigr| +\frac{\epsilon}{k}\Bigr) y +1 \Bigr\}.
\end{align*}
Since the discriminant of the right hand side is 
\[ 
\Bigl( \bigl| \frac{\cos(\pi\alpha)}{\cos(\pi\alpha/2)} \bigr| +\frac{\epsilon}{k}\Bigr)^2 - 
\Bigl( \frac{1}{\cos^2 (\pi\alpha /2)} - \frac{\epsilon}{k^2} \Bigr) = 
-4 \sin^2\frac{\pi\alpha}{2} + \frac{\epsilon}{k} \Bigl( 2 
\bigl| \frac{\cos(\pi\alpha)}{\cos(\pi\alpha/2)} \bigr| + \frac{1+\epsilon}{k} \Bigr), 
\] 
we begin with $\epsilon >0$ which makes this discriminant $<0$. 
Then, the absolute value of the integrand in \eqref{eq:3-2-7} is bounded by 
\[ 
e^{ka \bigl| \frac{\cos(\pi\alpha)}{\cos(\pi\alpha/2)} \bigr|} 2\sin\frac{\pi\alpha}{2}\, \frac{1}{k^2} 
\Bigl\{ \Bigl( \frac{1}{\cos^2 (\pi\alpha /2)} - \frac{\epsilon}{k^2} \Bigr) y^2 
-2 \Bigl( \bigl| \frac{\cos(\pi\alpha)}{\cos(\pi\alpha/2)} \bigr| +\frac{\epsilon}{k}\Bigr) y +1 \Bigr\}^{-1}, 
\] 
which is integrable in $y$ and independent of $n$. 
The poitwise limit of the integrand is seen from \eqref{eq:3-2-8}. 
Consequently, setting $\lim_{n\to\infty} \beta_n = \beta \in \{0, 1, \infty\}$, we have 
\begin{align}
\lim_{n\to\infty} f(k, n, t\theta_n) 
&= \frac{2\sin (\pi\alpha /2)}{\pi\alpha} \int_0^\infty e^{-t\beta (ky)^{\frac{1}{\alpha}}} 
\bigl\{ \frac{1}{\cos^2 (\pi\alpha /2)} y^2 + 2 \frac{\cos(\pi\alpha)}{\cos (\pi\alpha /2)} y +1 \bigr\}^{-1}dy 
\notag\\ 
&= \frac{\sin(\pi\alpha)}{\pi\alpha} \int_0^\infty \frac{e^{-t\beta (k\xi\cos(\pi\alpha/2))^{\frac{1}{\alpha}}}}
{\xi^2 + 2\xi \cos(\pi\alpha) +1} d\xi.
\label{eq:3-2-9}
\end{align}
Note that \eqref{eq:3-2-9} equals $1$ or $0$ according to $\beta =0$ or $\infty$.

Secondly we treat the case of $\alpha = 1/2$ in \eqref{eq:3-2-5}. 
Setting $\alpha = 1/2$ in \eqref{eq:3-2-3}, we have
\begin{align}
f(k,n,s) &= \frac{2k}{\pi n} \int_0^\infty \frac{e^{-sx^2}}{x} 
\frac{\sin(\sqrt{2}x)}{1-2(1-\frac{k}{n}) \cos(\sqrt{2}x) + (1-\frac{k}{n})^2} \,dx \notag\\ 
&= \frac{2k}{\pi} \int_0^\infty e^{-\frac{sy^2}{2n^2}} \,\frac{\sin (y/n)}{y/n}
\frac{1}{2n(n-k)(1-\cos (y/n)) +k^2} \, dy.
\label{eq:3-2-10}
\end{align}

\begin{lemma}\label{lem:3-3}
We have 
%
%\begin{equation}\label{eq:3-2-11}
\[ 
\frac{2k}{\pi} \int_0^{2\pi n} e^{-\frac{sy^2}{2n^2}} \,\frac{\sin (y/n)}{y/n}\, 
\frac{1}{2n(n-k)(1-\cos (y/n)) +k^2} dy < \;\text{\eqref{eq:3-2-10}}\; \leqq 1.
\] 
%\end{equation}
%
\end{lemma}
\textbf{Proof} \ 
The second inequality is obvious from \eqref{eq:2-7}, the definition of $f(k,n,s)$. 
For the first inequality, we divide the integral in \eqref{eq:3-2-10} as 
\[ 
\int_0^\infty = \sum_{j=1}^\infty \int_{2\pi nj}^{2\pi n(j+1)} \qquad 
\text{and note each} \qquad \int_{2\pi nj}^{2\pi n(j+1)} >0 \quad \text{for} \ j\in\mathbb{N}.
\] 
In fact, 
\begin{align*}
&\int_{2\pi nj}^{2\pi n(j+1)} e^{-\frac{sy^2}{2n^2}} \,\frac{\sin (y/n)}{y/n}\, 
\frac{1}{2n(n-k)(1-\cos (y/n)) +k^2}\, dy \\ 
&= \int_{2\pi j}^{2\pi (j+1)} e^{-\frac{sx^2}{2}}\, \frac{\sin x}{x}\,
\frac{1}{2(n-k)(1-\cos x) +k^2}\, dx 
= \int_{2\pi j}^{2\pi j+\pi} + \int_{2\pi j+\pi}^{2\pi j+2\pi} \\ 
&= \int_0^\pi\frac{\sin x}{2(n-k)(1-\cos x)+k^2} \Bigl\{ 
\frac{e^{-\frac{s}{2}(2\pi j+ x)^2}}{2\pi j +x} - \frac{e^{-\frac{s}{2}(2\pi j+2\pi - x)^2}}{2\pi j +2\pi -x}\Bigr\} dx
\end{align*}
where $\{~\cdot~\} >0$ for $0\leqq x < \pi$. 
\hfill $\blacksquare$

\bigskip

We compute the limit of 
\begin{equation}\label{eq:3-2-12}
f(k,n, t\theta_n) = \frac{2k}{\pi} \int_0^\infty e^{-t\beta_n y^2/2}\, \frac{\sin (y/n)}{y/n}\,
\frac{1}{2n(n-k)(1-\cos (y/n)) +k^2} dy
\end{equation}
($s =t\theta_n$, $\beta_n = \theta_n/n^2$). 
If $\beta_n \to \infty$ as $n\to\infty$, 
\[ 
\bigl| \text{\eqref{eq:3-2-12}} \bigr| \leqq \frac{2k}{\pi} \int_0^\infty e^{-t\beta_n y^2/2} 
\frac{1}{k^2}\, dy = \sqrt{\frac{2}{\pi}} \frac{1}{k\sqrt{t\beta_n}} \ \xrightarrow[n\to\infty] \ 0.
\] 
If $\beta_n\to 1$ as $n\to\infty$, since the integrand of \eqref{eq:3-2-12} is uniformly integrablein $n$, we have 
\[ 
\text{\eqref{eq:3-2-12}} \ \xrightarrow[n\to\infty] \ \frac{2k}{\pi} 
\int_0^\infty e^{-ty^2/2} \frac{1}{y^2+k^2} dy, 
\] 
which agrees with \eqref{eq:3-2-9} for $\alpha =1/2$ and $\beta =1$. 
Finally, let $\beta_n \to 0$ as $n\to\infty$. 
By Lemma~\ref{lem:3-3}, 
\[ 
\frac{2k}{\pi} \int_0^{2\pi n} e^{-t\beta_n y^2/2}\, \frac{\sin (y/n)}{y/n}\,
\frac{1}{2n(n-k)(1-\cos (y/n)) +k^2} dy < \text{\eqref{eq:3-2-12}} \leqq 1.
\] 
We show the above leftmost side tends to $1$ as $n\to\infty$. 
For any $\delta\in (0, \pi)$, 
\begin{align*}
&\Bigl| \int_{\delta n}^{2\pi n-\delta n} e^{-t\beta_ny^2/2} \frac{\sin(y/n)}{y/n} \frac{1}{2n(n-k)(1-\cos(y/n))+k^2} 
dy \Bigr| \\ 
&\leqq \int_\delta^{2\pi -\delta} \frac{ndx}{2n(n-k)(1-\cos\delta)+k^2} \ \xrightarrow[n\to\infty] \ 0, \\ 
&\Bigl| \int_{2\pi n- \delta n}^{2\pi n} e^{-t\beta_ny^2/2} \frac{\sin(y/n)}{y/n} \frac{1}{2n(n-k)(1-\cos(y/n))+k^2} 
dy \Bigr| \\ 
&\leqq \int_{2\pi-\delta}^{2\pi} \frac{-\sin x}{x} \frac{n\; dx}{2n(n-k)(1-\cos x) +k^2} \leqq 
\frac{n}{2\pi -\delta} \int_{-\delta}^0 \frac{-\sin x\; dx}{2n(n-k)(1-\cos x)+k^2} \\ 
&= \frac{n}{2\pi-\delta}\, \frac{1}{2n(n-k)} \bigl\{ \log \bigl( 2n(n-k)(1-\cos\delta)\bigr) - \log k^2\bigr\} 
\ \xrightarrow[n\to\infty] \ 0,
\end{align*}
and 
\begin{align*}
&\int_0^{\delta n} e^{-t\beta_ny^2/2} \frac{\sin(y/n)}{y/n} \frac{1}{2n(n-k)(1-\cos(y/n))+k^2} dy \\ 
&\ \xrightarrow[n\to\infty] \ \int_0^\infty \frac{1}{y^2+k^2} dy = \frac{k}{2\pi}
\end{align*}
since uniform integrability of the integrand follows, by taking $\delta >0$ small enough, from 
\[ 
0<y<\delta n \ \Longrightarrow \ 2n(n-k)(1-\cos\frac{y}{n}) +k^2 \geqq 
\frac{2}{3}\bigl( 1-\frac{k}{n}\bigr) y^2 +k^2.
\] 
This completes the proof of Proposition~\ref{prop:2-2}.

\subsection*{Remark}

In order to describe time evolution of the limit shape (= macroscopic profile) $\omega_t$, 
we presented that of its transition measure $\mathfrak{m}_{\omega_t}$ in this paper. 
This expression enables us to read out the $t$-dependence of $\omega_t$ by way of the Markov transform \eqref{eq:r-1}. 
Although it is often difficult to write down a concrete formula for $\omega_t$, we can, for example, 
appeal to numerical computation to follow the evolution of the shape. 
It is surely important to seek a partial differential equation for $\omega_t$ itself in addition to \eqref{eq:1-18.5} 
for the Stieltjes transform of $\mathfrak{m}_{\omega_t}$. 
Another promising way is given by the logarithmic energy. 
It is known that the limit shape $\varOmega$ of Vershik--Kerov and Logan--Shepp is the unique minimizer 
of the following functional on the continuous diagrams $\mathbb{D}$: 
\[ 
\Theta (\omega) = 1+ \frac{1}{2} \iint_{\{s>t\}} \bigl( 1-\omega^\prime(s)\bigr)\bigl( 1+\omega^\prime(t)\bigr)
\log (s-t) dsdt 
\] 
with $\Theta(\varOmega) =0$ (see \cite{Ke03} and also \cite{Ho16}). 
Since our limit shape $\omega_t$ converges to $\varOmega$ in $\mathbb{D}$ as $t\to\infty$, 
it is interesting to ask whether $\Theta(\omega_t)$ decreases as $t$ goes  by (maybe for sufficiently large $t$).

In this paper we focus on the limit shape evolution (law of large numbers) without mentioning fluctuation 
(central limit theorem) of the macroscopic profile. 
For a dynamical aspect of such fluctuation for Young diagrams, see \cite{Fu16}. 
As algebraic and systematical approach to static concentration and fluctuation problems for Young diagrams, we refer to 
\cite{IvOl02}, \cite{Sn06} and \cite{DoFe16}.

\subsection*{Appendix}

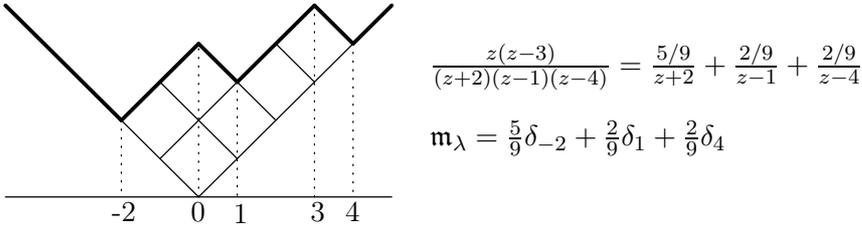
\begin{figure}[hbt]
{\unitlength 0.1in%
\begin{picture}(22.0000,10.4000)(4.0000,-14.4000)%
% LINE 2 0 3 0 Black White  
% 6 400 400 1400 1400 1400 1400 2400 400 2400 1400 400 1400
% 
\special{pn 8}%
\special{pa 400 400}%
\special{pa 1400 1400}%
\special{fp}%
\special{pa 1400 1400}%
\special{pa 2400 400}%
\special{fp}%
\special{pa 2400 1400}%
\special{pa 400 1400}%
\special{fp}%
% LINE 0 0 3 0 Black White  
% 12 2400 400 2200 600 2200 600 2000 400 2000 400 1600 800 1600 800 1400 600 1400 600 1000 1000 1000 1000 400 400
% 
\special{pn 20}%
\special{pa 2400 400}%
\special{pa 2200 600}%
\special{fp}%
\special{pa 2200 600}%
\special{pa 2000 400}%
\special{fp}%
\special{pa 2000 400}%
\special{pa 1600 800}%
\special{fp}%
\special{pa 1600 800}%
\special{pa 1400 600}%
\special{fp}%
\special{pa 1400 600}%
\special{pa 1000 1000}%
\special{fp}%
\special{pa 1000 1000}%
\special{pa 400 400}%
\special{fp}%
% LINE 2 0 3 0 Black White  
% 8 2000 800 1800 600 1800 1000 1600 800 1600 800 1200 1200 1600 1200 1200 800
% 
\special{pn 8}%
\special{pa 2000 800}%
\special{pa 1800 600}%
\special{fp}%
\special{pa 1800 1000}%
\special{pa 1600 800}%
\special{fp}%
\special{pa 1600 800}%
\special{pa 1200 1200}%
\special{fp}%
\special{pa 1600 1200}%
\special{pa 1200 800}%
\special{fp}%
% LINE 2 2 3 0 Black White  
% 10 1000 1000 1000 1400 1400 1400 1400 600 1600 800 1600 1400 2000 1400 2000 400 2200 600 2200 1400
% 
\special{pn 8}%
\special{pa 1000 1000}%
\special{pa 1000 1400}%
\special{dt 0.045}%
\special{pa 1400 1400}%
\special{pa 1400 600}%
\special{dt 0.045}%
\special{pa 1600 800}%
\special{pa 1600 1400}%
\special{dt 0.045}%
\special{pa 2000 1400}%
\special{pa 2000 400}%
\special{dt 0.045}%
\special{pa 2200 600}%
\special{pa 2200 1400}%
\special{dt 0.045}%
% STR 2 0 3 0 Black White  
% 4 2600 500 2600 600 1 0 0 0
% $\frac{z(z-3)}{(z+2)(z-1)(z-4)}=\frac{5/9}{z+2}+\frac{2/9}{z-1}+\frac{2/9}{z-4}$
\put(26.0000,-6.0000){\makebox(0,0)[lt]{$\frac{z(z-3)}{(z+2)(z-1)(z-4)}=\frac{5/9}{z+2}+\frac{2/9}{z-1}+\frac{2/9}{z-4}$}}%
% STR 2 0 3 0 Black White  
% 4 2600 900 2600 1000 1 0 0 0
% $\mathfrak{m}_\lambda=\frac{5}{9}\delta_{-2}+\frac{2}{9}\delta_1+\frac{2}{9}\delta_4$
\put(26.0000,-10.0000){\makebox(0,0)[lt]{$\mathfrak{m}_\lambda=\frac{5}{9}\delta_{-2}+\frac{2}{9}\delta_1+\frac{2}{9}\delta_4$}}%
% STR 2 0 3 0 Black White  
% 4 1360 1330 1360 1430 1 0 0 0
% 0
\put(13.6000,-14.3000){\makebox(0,0)[lt]{0}}%
% STR 2 0 3 0 Black White  
% 4 2160 1330 2160 1430 1 0 0 0
% 4
\put(21.6000,-14.3000){\makebox(0,0)[lt]{4}}%
% STR 2 0 3 0 Black White  
% 4 1580 1340 1580 1440 1 0 0 0
% 1
\put(15.8000,-14.4000){\makebox(0,0)[lt]{1}}%
% STR 2 0 3 0 Black White  
% 4 950 1330 950 1430 1 0 0 0
% -2
\put(9.5000,-14.3000){\makebox(0,0)[lt]{-2}}%
% STR 2 0 3 0 Black White  
% 4 1980 1330 1980 1430 1 0 0 0
% 3
\put(19.8000,-14.3000){\makebox(0,0)[lt]{3}}%
\end{picture}}%
\caption{profile and transition measure}
\label{fig:2}
\end{figure}

\paragraph{Profile and transition measure} 
A Young diagram $\lambda$ is displayed in the $xy$ coordinates plane, each box being a $\sqrt{2}\times\sqrt{2}$ square 
(Figure~\ref{fig:2}). 
Transition measure $\mathfrak{m}_\lambda$ of $\lambda$ is formed from the interlacing valley-peak coordinates 
$(x_1<y_1<\cdots <y_{r-1}<x_r)$ for $\lambda$ through the partial fraction expansion 
\[ 
\frac{(z-y_1)\cdots (z-y_{r-1})}{(z-x_1)\cdots (z-x_r)} = 
\sum_{i=1}^r \frac{\mathfrak{m}_\lambda(\{ x_i\})}{z-x_i} = 
\int_{\mathbb{R}} \frac{1}{z-x} \mathfrak{m}_\lambda (dx).
\] 
%
%\begin{figure}[hbt]
%%\begin{center}
%\include{msj2019spring-fig}
%%\vspace{-5mm}
%\caption{profile and transition measure}
%\label{fig:2}
%%\end{cemter}
%\end{figure}

\paragraph{Markov transform} 
The correspondence between a Young diagram (or its profile) and its transition measure is extended to a 
continuous diagram $\omega\in\mathbb{D}$ and its transition measure $\mathfrak{m}_\omega$ by 
\begin{equation}\label{eq:r-1}
\frac{1}{z} \exp\Bigl\{ \int_{\mathbb{R}} \frac{1}{x-z} \Bigl( \frac{\omega(x)-|x|}{2}\Bigr)^\prime dx \Bigr\} = 
\int_{\mathbb{R}} \frac{1}{z-x} \mathfrak{m}_\omega (dx), \qquad z\in\mathbb{C}^+.
\end{equation}

\paragraph{Free convolution and free compression} 
Let $(A, \phi)$ be a pair of unital $\ast$ algebra $A$ over $\mathbb{C}$ and state $\phi$ of $A$. 
For self-adjoint $a\in A$ and probability $\mu$ on $\mathbb{R}$, we write as $a\sim\mu$ if 
$\phi(a^n) = \int_{\mathbb{R}} x^n \mu(dx)$ for any $n$. 
If $a, b\in A$ are free and $a\sim\mu$, $b\sim\nu$, then $a+b\sim\mu\boxplus\nu$. 
The free convolution $\mu\boxplus\nu$ is uniquely determined for arbitrary compactly supported 
probabilities $\mu$ and $\nu$ on $\mathbb{R}$. 
If $q\in A$ is a self-adjoint projection such that $c = \phi(q)\neq 0$, $a\in A$ with $a\sim\mu$ 
(compactly supported), and $a, q$ are free, then probability $\nu$ on $\mathbb{R}$ is determined in such a way 
that $qaq\sim\nu$ in $(qAq, c^{-1}\phi\big|_{qAq})$, which is called the free compression of $\mu$ 
and denoted by $\mu_c$. 
For any compactly supported $\mu$ and $0<c\leqq 1$, $\mu_c$ is uniquely determined. 
The free convolution and free compression for compactly supported probabilities on $\mathbb{R}$ are 
characterized in terms of their free cumulants by 
\[ 
R_k(\mu\boxplus\nu) = R_k(\mu) + R_k(\nu), \qquad 
R_k(\mu_c) = c^{k-1} R_k(\mu), \qquad k\in\mathbb{N}.
\]

\subsection*{Acknowledgment} 
The author thanks Professor Takahiro Hasebe for comments and instructions.

\end{document}